\documentclass[10pt, DIV=17]{scrartcl}
\usepackage{amsmath,amsthm}
\usepackage{txfonts,lmodern,mathtools}
\usepackage{algorithmic}
\usepackage{subcaption,tikz}
\usepackage{hyperref}
\usepackage{siunitx}

\usetikzlibrary{decorations.markings}

\newcommand \email [1]{\href{mailto:#1}{#1}}

\theoremstyle{definition}
\newtheorem{alg}{Algorithm}[section]
\newtheorem{df}[alg]{Definition}

\theoremstyle{remark}
\newtheorem*{rmrk}{Remark}
\theoremstyle{plain}

\newcommand \dir      {{\mathrm D}} 
\newcommand \dirbound {{\Gamma_\dir}}
\newcommand \dirBC    {u_\dir}
\newcommand \neu      {{\mathrm N}}
\newcommand \neubound {{\Gamma_\neu}}
\newcommand \neunormal{\nu_\neu}
\newcommand \neuBC    {g}
\newcommand \RHS      {f}
\renewcommand \d      {\,\mathrm d}
\let \sei \coloneqq
\newcommand \nat      {\mathbb N}
\newcommand \reell    {\mathbb R}
\newcommand \cB       {\mathcal B}
\newcommand \cM       {\mathcal M}
\newcommand \cQ       {\mathcal Q}
\newcommand \MCL      {\mathbb M} 
\newcommand \BCL      {\mathbb B} 
\newcommand \T        {\mathsf T}
\newcommand \V        {\mathsf V}
\newcommand \Th        {\mathsf T\!_\mathsf h}
\newcommand \Tv        {\mathsf T\!_\mathsf v}
\newcommand \hsk        {\mathsf {hSk}}
\newcommand \vsk        {\mathsf {vSk}}
\newcommand \Q        {\textsl{\texttt Q}}
\newcommand \uni       {_{\textsc{uniform}}}
\newcommand \hb       {_{\textsc{hb}}}
\newcommand \thb      {_{\textsc{thb}}}
\newcommand \tsp      {_{\textsc{tspline}}}
\newcommand \p  {^+}
\newcommand \pp {^{+\!+}}
\newcommand \cN       {\mathcal N}
\newcommand \thbp     {_{\textsc{thb}\p}}
\newcommand \thbpp     {_{\textsc{thb}\pp}}

\newcommand \tsppp     {_{\textsc{tspline}\pp}}
\DeclareMathOperator \bisect {\textsc{bisect}}
\DeclareMathOperator \subdiv {\textsc{subdivide}}
\DeclareMathOperator \refine {\textsc{refine}}
\DeclareMathOperator \refuni {\textsc{refine\_uniform}}
\DeclareMathOperator \refthb {\textsc{refine\_thb}}
\DeclareMathOperator \closthb {\textsc{closure\_thb}}
\DeclareMathOperator \reftj {\textsc{refine\_tjunc}}
\DeclareMathOperator \reftel {\textsc{refine\_telem}}
\DeclareMathOperator \reftsp {\textsc{refine\_tspline}}
\DeclareMathOperator \clostsp {\textsc{closure\_tspline}}
\DeclareMathOperator \Span {span}
\DeclareMathOperator \conv {conv}
\DeclareMathOperator \supp {supp}
\DeclareMathOperator \ext {ext}
\newcommand \ceilfrac [2] {\bigl\lceil\tfrac{#1}{#2}\bigr\rceil}
\newcommand \floorfrac [2] {\bigl\lfloor\tfrac{#1}{#2}\bigr\rfloor}

\newcommand \tbigcup {\mathchoice{{\textstyle\bigcup}}{\bigcup}{\bigcup}{\bigcup}}

\newcommand \little {greedy }  
\newcommand \much {safe } 
\newcommand \Little {Greedy }  
\newcommand \Much {Safe }
\newcommand \leveldomain [2][\cQ] {\Omega_{#1,#2}}

\title{Adaptive Mesh Refinement Strategies in Isogeometric Analysis - A Computational Comparison%
\footnote{The authors gratefully acknowledge support by the Deutsche Forschungsgemeinschaft in the Priority Program 1748 "Reliable simulation techniques in solid mechanics: Development of non-standard discretization methods, mechanical and mathematical analysis" under the project "Adaptive isogeometric modeling of propagating strong discontinuities in heterogeneous materials" (KA3309/3-1 and PE2143/2-1).}}
\author{%
\renewcommand{\thefootnote}{\alph{footnote}}%
Paul Hennig\footnotemark[1], Markus K\"astner\footnotemark[1],\\
\renewcommand{\thefootnote}{\alph{footnote}}%
Philipp Morgenstern\footnotemark[2], Daniel Peterseim\footnotemark[2]%
}

\begin{document}
\maketitle
\renewcommand{\thefootnote}{\alph{footnote}}%
\footnotetext[1]{Technische Universit\"at Dresden, Institute for Solid Mechanics, \\ George-B\"ahr-Stra\ss{}e 3c, 01062 Dresden, Germany \\ email addresses: \email{paul.hennig@tu-dresden.de}, \email{markus.kaestner@tu-dresden.de}}
\footnotetext[2]{Rheinische Friedrich-Wilhelms-Universit\"at Bonn, Institute for Numerical Simulation, \\ Wegelerstr.\ 6, 53115 Bonn, Germany \\ email addresses: \email{morgenstern@ins.uni-bonn.de}, \email{peterseim@ins.uni-bonn.de}}
\begin{abstract}
We explain four variants of an adaptive finite element method with cubic splines and compare their performance in simple elliptic model problems. The methods in comparison are Truncated Hierarchical B-splines with two different refinement strategies, T-splines with the refinement strategy introduced by Scott et al.\ in 2012, and T-splines with an alternative refinement strategy introduced by some of the authors. In four examples, including singular and non-singular problems of linear elasticity and the Poisson problem, the H1-errors of the discrete solutions, the number of degrees of freedom as well as sparsity patterns and condition numbers of the discretized problem are compared.
\end{abstract}

\section{Introduction}

Adaptive Isogeometric Methods (AIGM) have gained widespread interest and are a very active field of research, investigating a wide range of refinement strategies. The ``usual'' mesh refinement, entitled $h$-refinement, currently competes with $p$-refinement (augmenting the polynomial degree), $k$-refinement (a particular combination of $h$- and $p$-refinement), $r$-refinement (redesigning the mesh) and their combinations. However, if B-splines or NURBS 
are considered as a basis, their tensor product nature will prohibit a truly local $h$-refinement within a single patch, and various approaches have been developed to overcome this restrictive tensor product structure.

The concept of T-splines as an $h$-refinement technique caught much attention \cite{bazilevs2010, doerfel2010}, but also incorporated algorithmic difficulties \cite{buffa2010}, in particular linear dependencies between the T-spline functions that should serve as a spline basis, and an unclear nesting behaviour of the generated spline spaces. Most of these problems could be solved in the last years \cite{daveiga2012, sederberg2012, scott2012, scott2014, mp2015, morgenstern2015}. 

Hierarchical (H)B-splines, a further promising $h$-refinement technique, have been introduced already in 1988 \cite{forsey1988} and developed to meet the requirements of isogeometric analysis \cite{kraft1997,Vuong2011,kvvb2014}. They have been enhanced to truncated hierarchical (TH)B-splines \cite{giannelli2012} by reducing the interaction between basis functions of different refinement levels. As a result, THB-splines improve the conditioning and reduce the bandwidth as well as the total number of non-zero entries in the system matrices \cite{Johannessen2015,Schillinger2012}. The mentioned interaction can be further regulated by the introduction of $m$-admissible meshes, where interacting basis functions in an element belong to at most $m$ different levels \cite{buffa2015a}.

Further $h$-refinement techniques are locally refined (LR-)B-splines \cite{dokken2013,bressan2013}, hierarchical T-splines \cite{evans2015}, modified T-splines \cite{kang2013}, PHT-splines \cite{deng2008,wang2011c} amongst many others. 
\begin{figure}
 \centering
 \includegraphics[scale=0.65]{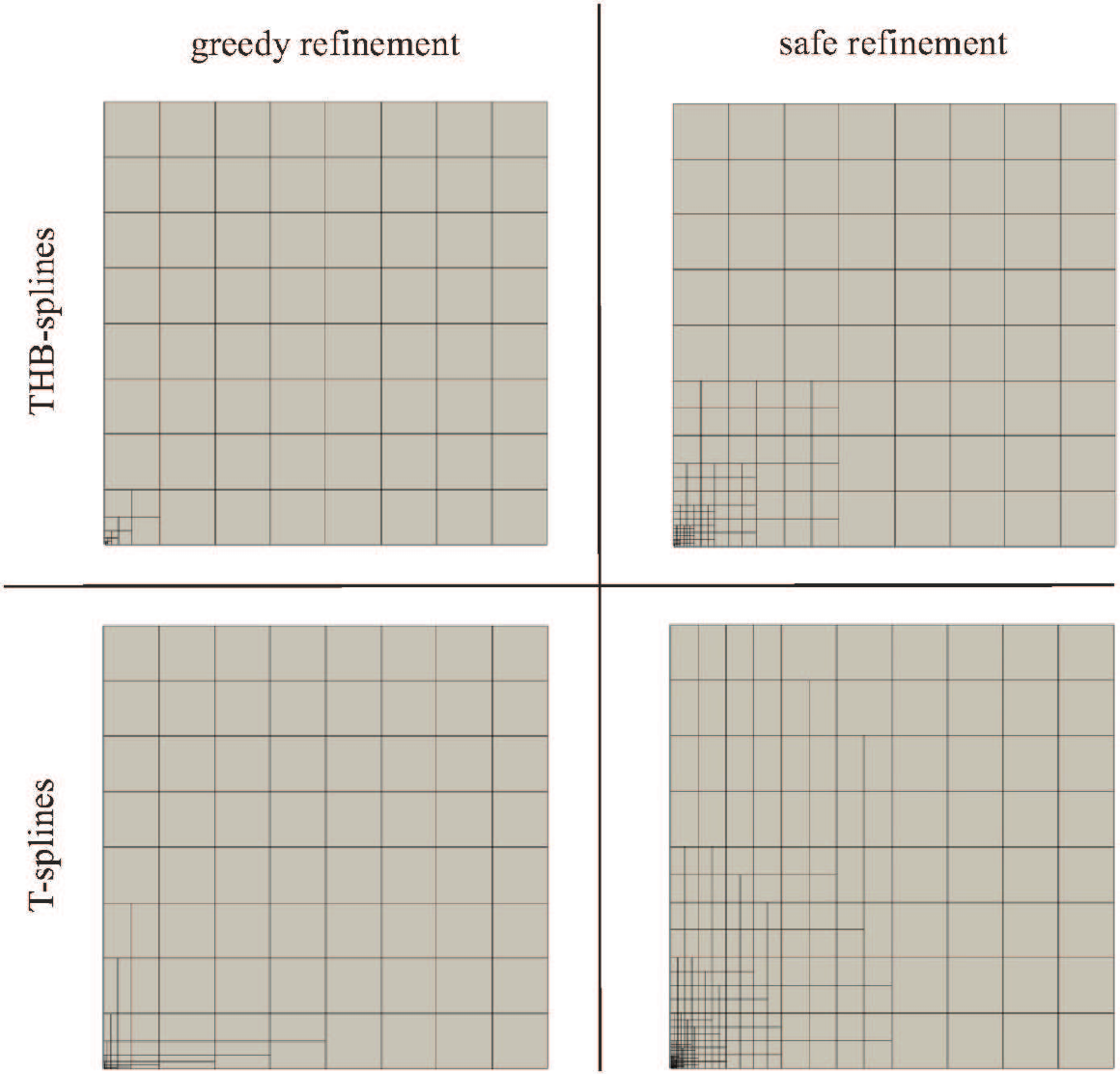}
 \caption{Refinement strategies: An initial square mesh with 64 elements is locally refined in the lower left corner using THB- and T-splines with a \little and a \much refinement. The illustrated meshes are the B\'ezier meshes (see Definition~\ref{df: bezier mesh}).}
 \label{fig:SquareCorner_Mesh}
\end{figure}

In this contribution, four different realizations of an Adaptive Isogeometric Method with $h$-refinement are compared: 
\begin{enumerate}
\item A refinement based on T-splines \cite{scott2012}, where the refinement process is divided into two steps. At first, marked elements are refined, and secondly, an additional refinement is processed to recover the linear independence of the T-spline functions, called \emph{analysis-suitability}.
\item A refinement based on T-splines \cite{morgenstern2015,mp2015}, where also the vicinity of the marked elements is considered. By defining a class of admissible T-meshes, the proposed refinement preserves the analysis-suitability of the T-splines directly.
\item A refinement based on THB-splines, where the mesh allows only a one-level difference between neighbouring mesh elements.
\item A refinement based on THB-splines, where only 2-admissible meshes are allowed.
\end{enumerate}
We refer to these methods as \emph{\little refinement} (method 1 and 3) and \emph{\much refinement} (method 2 and 4).

Some of these methods allow for a mathematical proof of linear complexity \cite{bgmp2016,mp2015}. Together with results on the convergence of the Adaptive Algorithm \cite{buffa2015a}, this allows for a proof of optimal convergence rates \cite{carstensen2014}.
This paper will compare the refinement of T-splines and THB-splines with a focus on the influence of the different mesh classes on the numerical solution and properties. In four examples, including singular and non-singular problems of linear elasticity and the Poisson problem, the H1-errors of the discrete solutions, the number of degrees of freedom as well as sparsity patterns and condition numbers of the discretized problem are compared.

To enable a unified implementation of the different refinement techniques, B\' ezier extraction is used. B\' ezier extraction provides a canonical form to use isogeometric analysis with different spline bases, has a strict element viewpoint, allows for an implementation into existing finite element codes and has been already developed for T-splines \cite{Borden2011} and THB-splines \cite{Hennig2016}.   

The paper is organized as follows. Section~\ref{sec: mesh refinement} introduces the refinement strategies to be compared. Section~\ref{sec: model problem} describes the problems that will be solved numerically using each of the presented methods. Section~\ref{sec: adaptive loop} investigates the adaptive algorithm and summarizes background theory on optimal convergence rates. The computational comparison is performed in Section~\ref{sec: comparison}, and conclusions are given in Section~\ref{sec: conclusion}.

\section{Mesh Refinement Strategies}\label{sec: mesh refinement}
In this section, we define several $h$-refinement strategies for bivariate B-splines. They will be described in the style 
\[\refine(\cQ,\cM) = (\tilde\cQ,\tilde\cB),\]
with $\cQ$ being a rectangular mesh in the index domain 
and $\cM\subseteq\cQ$ a set of elements (rectangles) to be refined. $\tilde\cQ$ and $\tilde\cB$ are the new refined mesh and the set of B-spline basis functions associated to that new mesh, respectively.

\subsection{Uniform refinement} 
\begin{figure}[ht]
\captionsetup[subfigure]{labelformat=empty}
\centering
\begin{subfigure}{.2\textwidth}
\centering
\begin{tikzpicture}
\draw (0,0) grid (3,2);
\end{tikzpicture}
\caption{$\cQ_0$}
\end{subfigure}\hspace{.05\textwidth}
\begin{subfigure}{.2\textwidth}
\centering
\begin{tikzpicture}[scale=.5]
\draw (0,0) grid (6,4);
\end{tikzpicture}
\caption{$\cQ_1$}
\end{subfigure}\hspace{.05\textwidth}
\begin{subfigure}{.2\textwidth}
\centering
\begin{tikzpicture}[scale=.25]
\draw (0,0) grid (12,8);
\end{tikzpicture}
\caption{$\cQ_2$}
\end{subfigure}\hspace{.05\textwidth}
\begin{subfigure}{.2\textwidth}
\centering
\begin{tikzpicture}[scale=.125]
\draw (0,0) grid (24,16);
\end{tikzpicture}
\caption{$\cQ_3$}
\end{subfigure}
\caption{Example for the uniform meshes $\cQ_0,\dots,\cQ_3$ for $M=3$ and $N=2$.}
\label{fig: uniform meshes}
\end{figure}
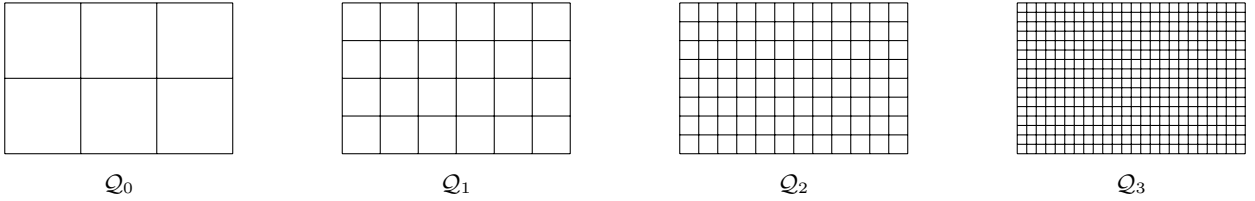

We assume the initial mesh $\cQ_0$ to be a tensor product mesh, and its elements are closed squares with side length~1 (see Figure~\ref{fig: uniform meshes}),
\[\cQ_0\sei\Bigl\{[m-1,m]\times[n-1,n]\mid m\in\{1,\dots,M\},n\in\{1,\dots,N\}\Bigr\}.\]
The corresponding spline basis $\cB_0$ is spanned by the corresponding tensor-product B-splines.
For each level $k\in\nat_0$, we define the tensor-product mesh
\[\cQ_k\sei\Bigl\{[x-2^{-k},x]\times[y-2^{-k},y]\mid 2^k x\in\{1,\dots,2^k M\},2^k y\in\{1,\dots,2^k N\}\Bigr\}\]
and the corresponding spline space $\cB_k$ of tensor-product bivariate B-spline basis functions.

\begin{df}
We define the uniform refinement routine by
\[\refuni(\cQ_k,\cM) \sei (\cQ_{k+1},\cB_{k+1})\quad\text{for any }k\in\nat_0\text{ and }\cM\subseteq\cQ.\]
Note that the set of marked elements $\cM$ enters only for formal reasons and has no effect on the refinement.
We denote the class of uniform meshes by $\MCL\uni\sei\{\cQ_n\mid n\in\nat_0\}$.
\end{df}
\begin{df}[level]
Given $k\in\mathbb N$ and $\Q\in\cQ_k$, we denote the \emph{level} of $\Q$ by $\ell(\Q)\sei k$.
\end{df}

\subsection{Truncated Hierarchical B-splines} \label{sec:Truncated Hierarchical B-splines}
\begin{figure}
\captionsetup[subfigure]{labelformat=empty}
\centering
\begin{subfigure}{.25\textwidth}
\includegraphics[width=\textwidth]{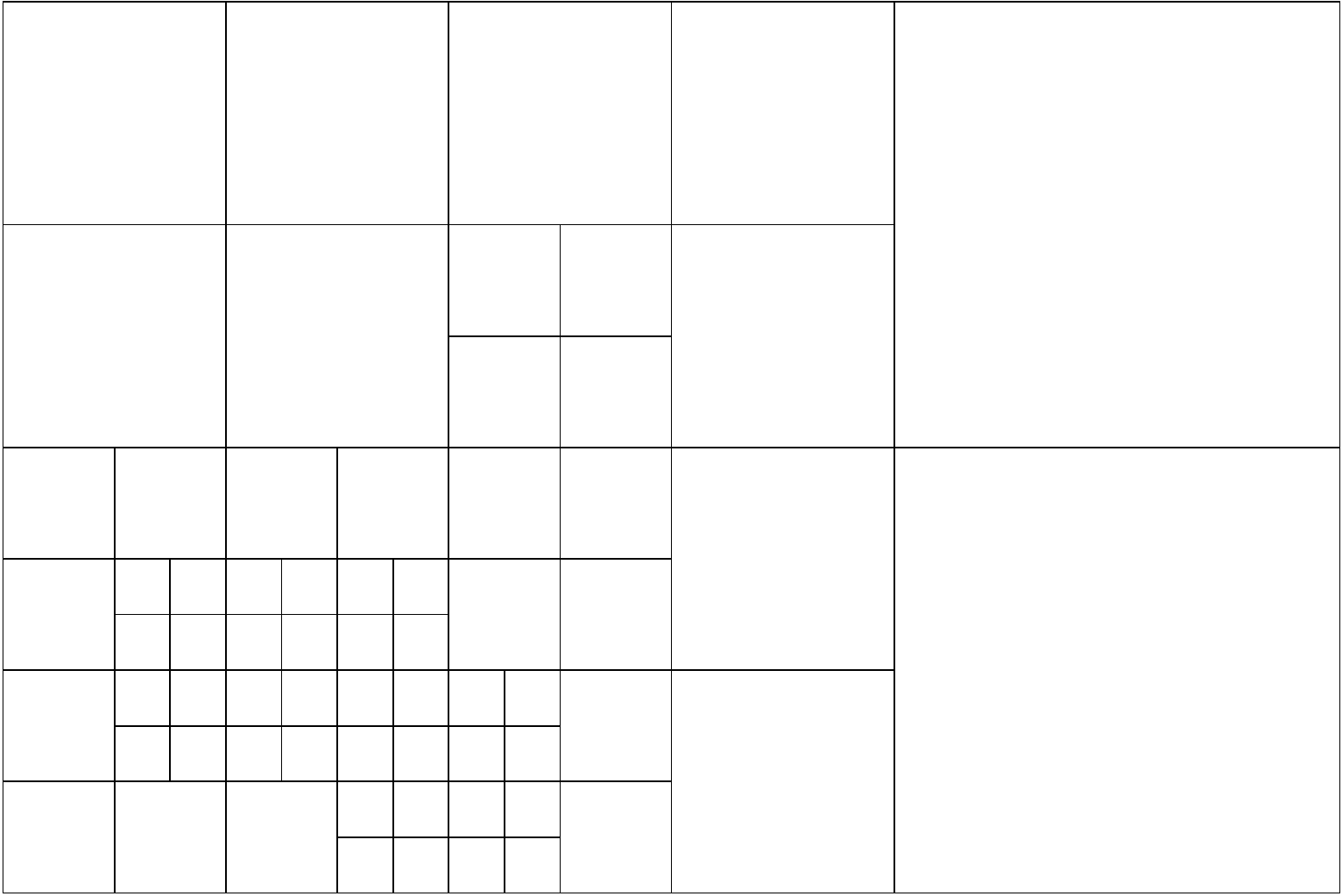}
\caption{$\cQ$}
\end{subfigure}\\[1ex]
\begin{subfigure}{.2\textwidth}
\includegraphics[width=\textwidth]{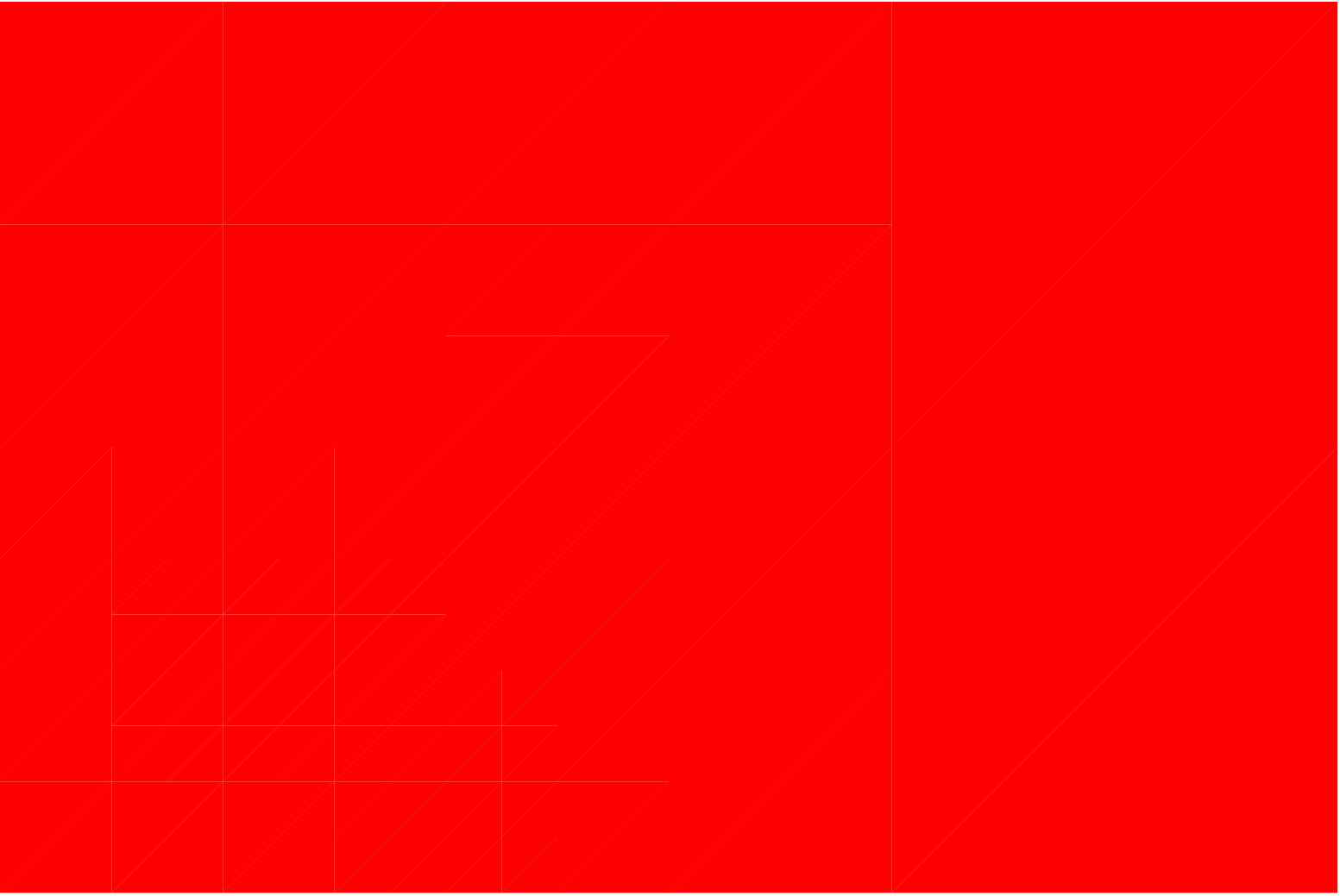}
\caption{$\leveldomain0$}
\end{subfigure}\hspace{.05\textwidth}
\begin{subfigure}{.2\textwidth}
\includegraphics[width=\textwidth]{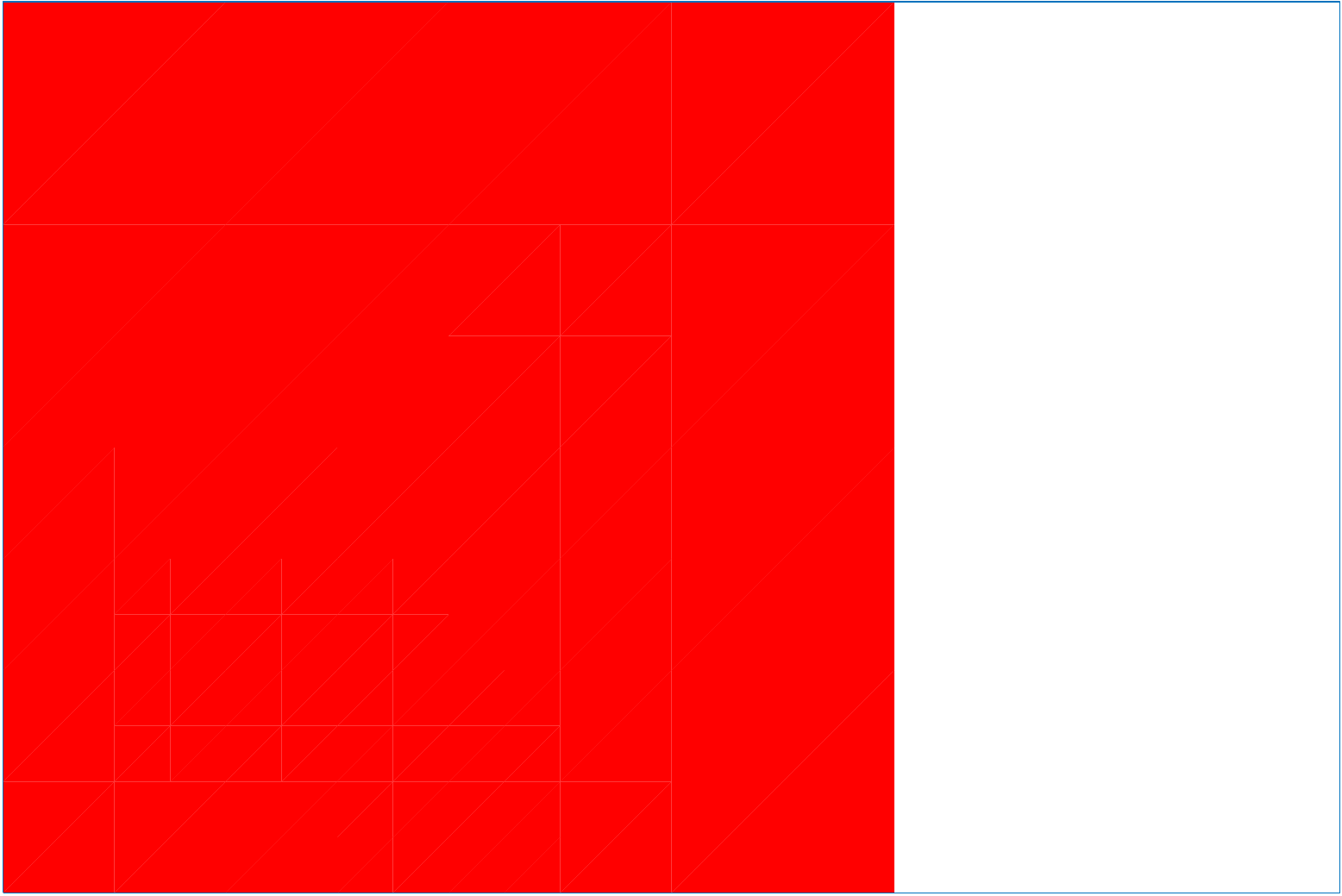}
\caption{$\leveldomain1$}
\end{subfigure}\hspace{.05\textwidth}
\begin{subfigure}{.2\textwidth}
\includegraphics[width=\textwidth]{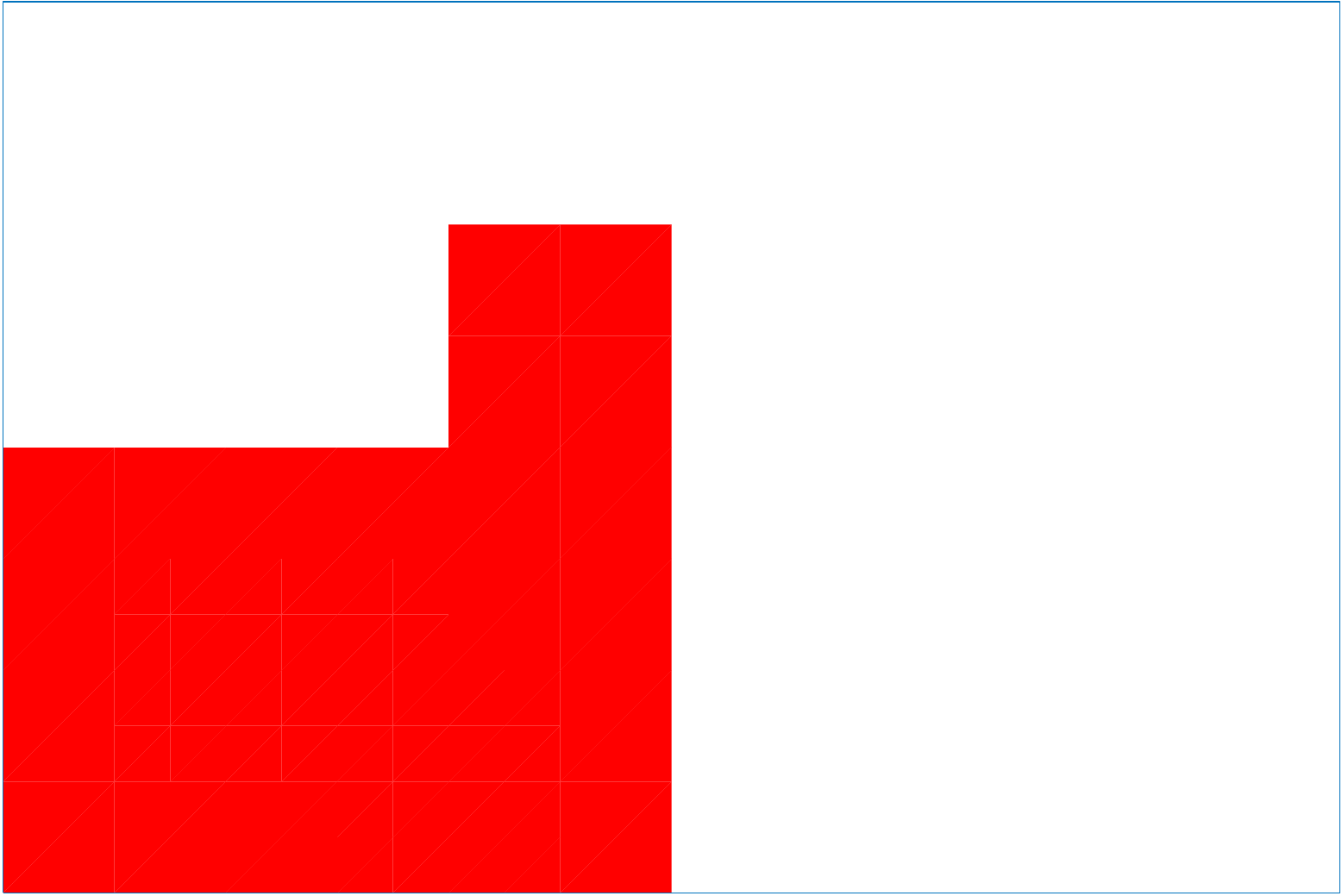}
\caption{$\leveldomain2$}
\end{subfigure}\hspace{.05\textwidth}
\begin{subfigure}{.2\textwidth}
\includegraphics[width=\textwidth]{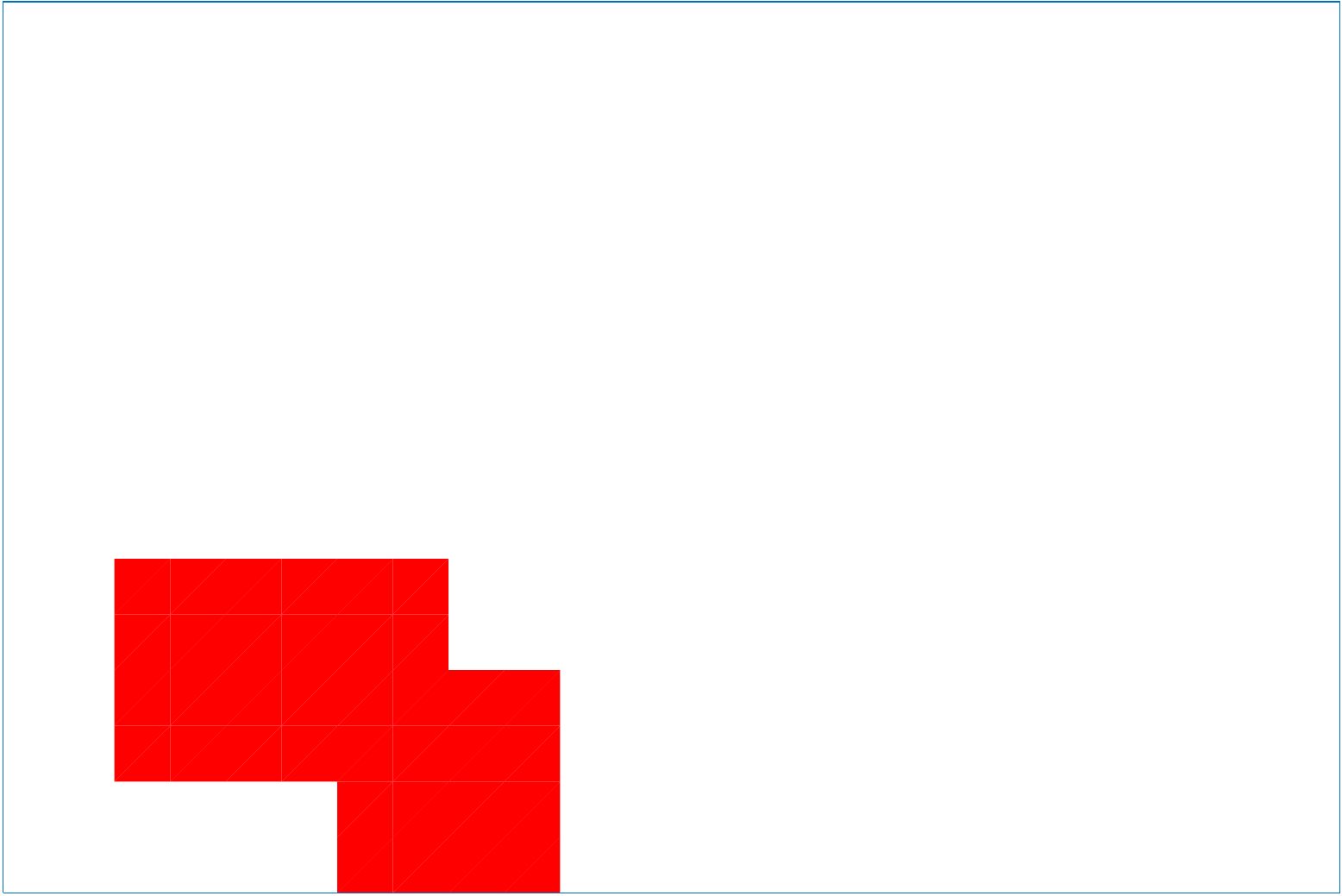}
\caption{$\leveldomain3$}
\end{subfigure}
\caption{Example for level-$k$ domains, for $k=0,\dots,3$. The domains $\leveldomain0,\dots,\leveldomain3$ are shaded in red.}
\label{fig: level-k domains}
\end{figure}

For the use of Truncated Hierarchical B-splines (THB-splines), the underlying rectangular mesh $\cQ$ may consist of finitely many elements from meshes in $\MCL\uni$, such that any two elements of $\cQ$ have disjoint interior, and the union of all elements of $\cQ$ is the same domain $[0,M]\times[0,N]$ that is covered by uniform meshes. In particular, $\cQ$ is meant to contain elements of different levels:
\[\MCL_\mathrm{THB}\sei\Bigl\{\cQ\subset\tbigcup_{\cQ'\in\MCL\uni}\cQ'\mid\#\cQ<\infty,\ \tbigcup\cQ=[0,M]\times[0,N],\ \forall \Q,\Q'\in\cQ:\ \operatorname{int}(\Q)\cap\operatorname{int}(\Q')=\emptyset\Bigr\}.\]
\begin{df}[level-$k$ domain]
Given some mesh $\cQ\in\MCL_\mathrm{THB}$, and $k\in\nat_0$, we denote by $\leveldomain k$ the domain that is covered by ``level-$k$ or finer'' elements, $\leveldomain k\sei\tbigcup\bigl\{\Q\in\cQ\mid\ell(\Q)\ge k\bigr\}$. See Figure~\ref{fig: level-k domains} for an example.
\end{df}

The classical (non-truncated) Hierarchical B-spline basis reads
\[\cB\hb^\cQ\sei\bigcup_{k\in\nat_0}\left\{B\in\cB_k\mid\supp B\subseteq\leveldomain k\text{ and }\supp B\nsubseteq\leveldomain{k+1}\right\}.\]
THB-splines involve an alternative choice of basis functions that span the same space as the basis $\cB\hb^\cQ$ above. These basis functions have reduced overlap compared to $\cB\hb^\cQ$ and hence provide sparser Galerkin and Collocation matrices when used for solving a discretized PDE. 
\begin{df}[Truncation of the classical basis, \cite{giannelli2012}]
Let $B=\sum_{\bar B\in\cB_j}c_{\bar B,B}\bar B\in\operatorname{span}\cB_j$, then 
\[\operatorname{trunc}^j_\cQ(B)\sei\sum_{\bar B\in\cB_{j}\setminus\cB\hb^\cQ}c_{\bar B,B}\bar B.\]
In addition, for any $B\in\cB_k$, $k\in\nat_0$ we define the successive truncation w.r.t.\ all higher levels
\[\operatorname{Trunc}_\cQ(B)\sei
\operatorname{trunc}_\cQ^K(\dots\operatorname{trunc}_\cQ^{k+2}(\operatorname{trunc}_\cQ^{k+1}(B))\dots)\]
with $K=\max\{k\in\nat_0\mid \cB_k\cap\cB\hb^\cQ\neq\emptyset\}$.
\end{df}

\begin{df}[Subdivision]
For $k\in\nat_0$ and $\Q\in\cQ_k$, we define
\[\subdiv(\Q)\sei \bigl\{\Q'\in\cQ_{k+1}\mid \Q'\subset \Q\bigr\}\]
and for $\cM\subset\cQ\in\MCL\thb$, we denote $\subdiv(\cM)\sei\bigcup_{\Q\in\cM}\subdiv(\Q)$.
\end{df}
\begin{df}
The refinement procedure for THB-splines reads
$\refthb(\cQ,\cM) \sei \bigl(\tilde\cQ,\cB\thb^{\tilde\cQ}\bigr)$, with
\begin{align*}
\tilde\cQ&\sei\cQ\setminus\cM\cup\subdiv(\cM)\\
\cB\thb^{\tilde\cQ}&\sei\left\{\operatorname{Trunc}(B)\mid B\in\cB\hb^{\tilde\cQ}\right\}.
\end{align*}
\end{df}
\subsubsection{\Little refinement for THB-Splines}
\begin{figure}[ht]
\centering
\begin{tikzpicture}[baseline=.875cm]
\fill[blue!50] (0,0) rectangle (.25,.25);
\draw (0,0) grid (2,2);
\draw[step=.5] (0,0) grid (1,1);
\draw[step=.25] (0,0) grid (.5,.5);
\end{tikzpicture}
\quad$\rightarrow$\quad
\begin{tikzpicture}[baseline=.875cm]
\fill[blue!50] (0,0) rectangle (.25,.25);
\draw (0,0) grid (2,2);
\draw[step=.5] (0,0) grid (1,1);
\draw[step=.25] (0,0) grid (.5,.5);
\end{tikzpicture}
\quad$\rightarrow$\quad
\begin{tikzpicture}[baseline=.875cm]
\draw (0,0) grid (2,2);
\draw[step=.5] (0,0) grid (1,1);
\draw[step=.25] (0,0) grid (.5,.5);
\draw[step=.125] (0,0) grid (.25,.25);
\end{tikzpicture}
\caption{Another example for the \little THB-spline refinement. First, an element $\Q\in\cQ$ is marked (highlighted in blue), hence $\cM=\{\Q\}$. Second, $\closthb\p(\cM,\cQ)$ is computed, which now coincides with the actually marked element. Third, all elements in $\closthb\p(\cM,\cQ)$ are subdivided, which is only $\Q$.}
\label{fig: ex 2 little thb refinement}
\end{figure}
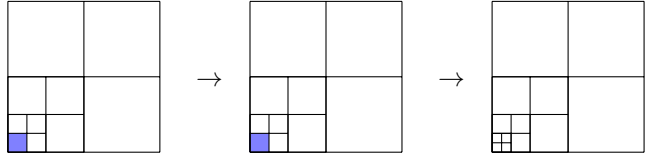

The \little THB-spline refinement defined below allows only a one-level difference between neighbouring mesh elements to produce graded meshes. Examples are given in Figure~\ref{fig: ex 1 little thb refinement} and Figure~\ref{fig: ex 2 little thb refinement}.
\begin{df}[\Little refinement for THB-Splines]
We define for each $\Q\in\cQ$ the \emph{coarse neighbourhood}
\begin{align*}
\cN\thbp(\Q)&\sei \left\{\Q'\in\cQ\mid \ell(\Q')>\ell(\Q),\ \Q\cap \Q'\neq\emptyset\right\}
\intertext{with generalized notations $\cN\thbp(\cM)\sei\bigcup_{\Q\in\cM}\cN\thbp(\Q)$ and $\cN\thbp^k(\cM)\sei
\smash{\underbrace{\cN\thbp(\dots\cN\thbp(}_{k\text{ times}}}\cM)\dots)$.
We further define the \emph{closure}}
\closthb\p(\cM,\cQ)&\sei \smash{ \bigcup_{k=0}^{\max\ell(\cM)}\cN\thbp^k(\cM), }
\intertext{and the extended refinement procedure}
\refthb\p(\cQ,\cM) &\sei \refthb\bigl(\cQ,\closthb\p(\cM,\cQ)\bigr).
\end{align*}
\end{df}
\begin{figure}[ht]
\centering
\begin{tikzpicture}[baseline=.675cm, scale=.75]
\fill[blue!50] (3,2.5) rectangle (3.25,2.75);
\draw (-1,-1) grid (4,3);
\draw[step=.5] (2,2) grid (4,3);
\draw[step=.25] (3,2.5) grid (3.5,3);
\end{tikzpicture}
\quad$\rightarrow$\quad
\begin{tikzpicture}[baseline=.675cm, scale=.75]
\fill[blue!50] (2,1)-|(4,2)-|(3.5,2.5)-|(3.25,2.75)-|(3,3)-|(2.5,2)-|(2,1);
\draw (-1,-1) grid (4,3);
\draw[step=.5] (2,2) grid (4,3);
\draw[step=.25] (3,2.5) grid (3.5,3);
\end{tikzpicture}
\quad$\rightarrow$\quad
\begin{tikzpicture}[baseline=.675cm, scale=.75]
\draw (-1,-1) grid (4,3);
\draw[step=.5] (2,1) grid (4,3);
\draw[step=.25] (2.5,2) grid (3.5,3);
\draw[step=.125] (3,2.5) grid (3.25,2.75);
\end{tikzpicture}
\caption{Example for the \little THB-spline refinement. First, an element $\Q\in\cQ$ is marked (highlighted in blue), hence $\cM=\{\Q\}$. Second, $\closthb\p(\cM,\cQ)$ is computed (highlighted in blue). Third, all elements in $\closthb\p(\cM,\cQ)$ are subdivided.}
\label{fig: ex 1 little thb refinement}
\end{figure}
\subsubsection{\Much refinement for THB-Splines}
The \much refinement for THB-splines defined below is conceptionally similar to the refinement procedure described in \cite{buffa2015a} and \cite{bgmp2016}. It only differs in the construction of the neighbourhood $\cN\thbpp$, where a different level is chosen for the B-splines whose supports are used in the construction. The resulting meshes are 2-admissible, meaning that interacting basis functions in an element belong to at most two different levels. Examples are given in Figure~\ref{fig: ex 1 much thb refinement} and Figure~\ref{fig: ex 2 much thb refinement}.
\begin{df}[\Much refinement for THB-Splines]
We define for each $\Q\in\cQ$ 
the \emph{same-level neighbourhood}
\begin{align*}
\cN_{\textsc{sl}}(\Q)&\sei \left\{\Q'\in\cQ_{\ell(\Q)}\mid \exists B\in\cB_{\ell(\Q)}:\ \Q\subset\supp B\supset \Q'\right\},
\shortintertext{the \emph{coarse neighbourhood}}
\cN\thbpp(\Q)&\sei \left\{\Q'\in\cQ\mid \ell(\Q')<\ell(\Q),\ \exists \Q''\in\cN_{\textsc{sl}}(\Q):\ \Q''\subset\Q' \right\},
\shortintertext{the \emph{closure}}
\closthb\pp(\cM,\cQ)&\sei \smash{ \bigcup_{k=0}^{\max\ell(\cM)}\cN\thbpp^k(\cM), }
\intertext{and the extended refinement procedure}
\refthb\pp(\cQ,\cM) &\sei \refthb\bigl(\cQ,\closthb\pp(\cM,\cQ)\bigr).
\end{align*}
\end{df}
\begin{figure}[ht]
\centering
\begin{tikzpicture}[baseline=.675cm, scale=.75]
\fill[blue!50] (3,2.5) rectangle (3.25,2.75);
\draw (-1,-1) grid (4,3);
\draw[step=.5] (1,1) grid (4,3);
\draw[step=.25] (3,2.5) grid (3.5,3);
\end{tikzpicture}
\quad$\rightarrow$\quad
\begin{tikzpicture}[baseline=.675cm, scale=.75]
\fill[blue!50] (0,0) rectangle (4,3);
\fill[white] (1,1) rectangle (4,3);
\fill[blue!50] (2,1.5) rectangle (4,3);
\fill[white] (3,2.5) rectangle (3.5,3);
\fill[blue!50] (3,2.5) rectangle (3.25,2.75);
\draw (-1,-1) grid (4,3);
\draw[step=.5] (1,1) grid (4,3);
\draw[step=.25] (3,2.5) grid (3.5,3);
\end{tikzpicture}
\quad$\rightarrow$\quad
\begin{tikzpicture}[baseline=.675cm, scale=.75]
\draw (-1,-1) grid (4,3);
\draw[step=.5] (0,0) grid (4,3);
\draw[step=.25] (2,1.5) grid (4,3);
\draw[step=.125] (3,2.5) grid (3.25,2.75);
\end{tikzpicture}
\caption{Example for the \much THB-spline refinement. First, an element $\Q\in\cQ$ is marked (highlighted in blue), hence $\cM=\{\Q\}$. Second, $\closthb\p(\cM,\cQ)$ is computed (highlighted in blue). Third, all elements in $\closthb\p(\cM,\cQ)$ are subdivided.}
\label{fig: ex 1 much thb refinement}
\end{figure}
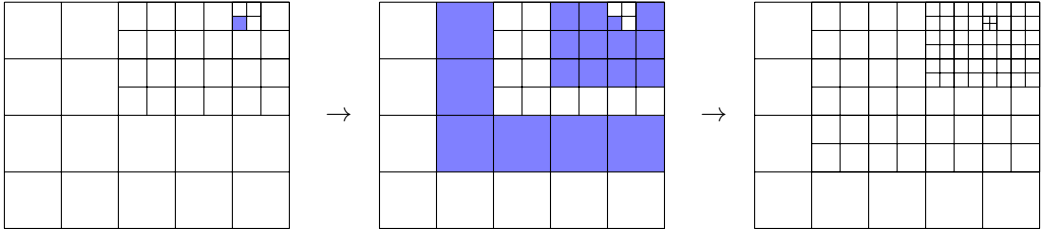
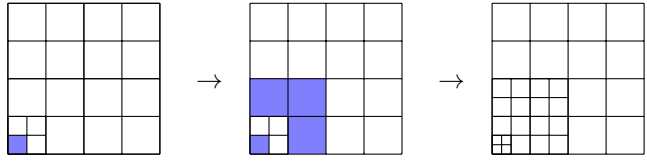
\begin{figure}[ht]
\centering
\begin{tikzpicture}[baseline=.875cm]
\fill[blue!50] (0,0) rectangle (.25,.25);
\draw (0,0) grid (2,2);
\draw[step=.5] (0,0) grid (2,2);
\draw[step=.25] (0,0) grid (.5,.5);
\end{tikzpicture}
\quad$\rightarrow$\quad
\begin{tikzpicture}[baseline=.875cm]
\fill[blue!50] (0,0) rectangle (1,1);
\fill[white] (0,0) rectangle (.5,.5);
\fill[blue!50] (0,0) rectangle (.25,.25);
\draw[step=.5] (0,0) grid (2,2);
\draw[step=.25] (0,0) grid (.5,.5);
\end{tikzpicture}
\quad$\rightarrow$\quad
\begin{tikzpicture}[baseline=.875cm]
\draw[step=.5] (0,0) grid (2,2);
\draw[step=.25] (0,0) grid (1,1);
\draw[step=.125] (0,0) grid (.25,.25);
\end{tikzpicture}
\caption{Another example for the \much THB-spline refinement. First, an element $\Q\in\cQ$ is marked (highlighted in blue), hence $\cM=\{\Q\}$. Second, $\closthb\p(\cM,\cQ)$ is computed. Third, all elements in $\closthb\p(\cM,\cQ)$ are subdivided.}
\label{fig: ex 2 much thb refinement}
\end{figure}

\subsection{T-splines}
While the refinement strategies for THB-splines presented above differ only in the choice of the neighbourhoods $\cN\thbp$ and $\cN\thbpp$, the refinement strategies for T-splines below are conceptionally different.
Throughout this paper, we denote the refinement procedure introduced in \cite{scott2012} as \emph{\little refinement for T-splines}.
It relies on T-junctions and T-junction extensions, and the set of children from a single element's bisection may be any set of two rectangles with disjoint interior such that their union is the parent element. This is, any element can be bisected in \emph{both} $x$- or $y$-direction, and the children may differ in size from each other.
On the other hand, the refinement strategy from \cite{mp2015}, denoted \emph{\much refinement for T-splines}, follows the structure of the THB-spline refinement above, marking coarser elements in the neighbourhood of marked elements, and then refining all marked elements at the same time. For the \much refinement, each element can be bisected in \emph{either} $x$- or $y$-direction, producing two rectangles with equal size, i.e., each element has a unique fixed set of children.

For the sake of legibility, we only give a definition of odd-degree T-splines. 
However, both refinement procedures $\reftsp\p$ and $\reftsp\pp$ are also suitable for even- or mixed-degree T-splines \cite{scott2012,mp2015,daveiga2013}.
\subsubsection{\Little refinement for T-splines}
\begin{df}
For any rectangle $\Q=[x,x+\tilde x]\times[y,y+\tilde y]$ and parameters $j\in\{1,2\}$, $0<q<1$, we define the refinement \[\bisect_{j,q}(\Q)\sei\begin{cases} 
\bigl\{ [x,x+q\tilde x]\times[y,y+\tilde y],\ [x+q\tilde x,x+\tilde x]\times[y,y+\tilde y]\bigr\}&\text{if }j=1,\\[.25ex]
\bigl\{ [x,x+\tilde x]\times[y,y+q\tilde y],\ [x,x+\tilde x]\times[y+q\tilde y,y+\tilde y]\bigr\}&\text{if }j=2.
\end{cases}\]
\end{df}
\begin{figure}[ht]
\centering
\begin{tikzpicture}[baseline=.4cm, scale=1.15]
\draw (0,0) rectangle (1,1);
\end{tikzpicture}
\quad$\xrightarrow{\bisect_{2,\,0.3}}$\quad
\begin{tikzpicture}[baseline=.4cm, scale=1.15]
\draw (0,0) rectangle (1,1) (0,.3)--(1,.3);
\end{tikzpicture}
\quad$\xrightarrow{\bisect_{1,\,0.6}}$\quad
\begin{tikzpicture}[baseline=.4cm, scale=1.15]
\draw (0,0) rectangle (1,1) (0,.3)--(1,.3) (.6,.3)--(.6,1);
\end{tikzpicture}
\caption{Example for $\bisect_{j,q}$.}
\end{figure}
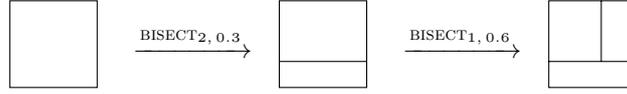
\begin{df}
We define the mesh class $\MCL\tsp$ inductively through the $\bisect_{j,q}$ routine;
\begin{align*}
\cQ_0&\in\MCL\tsp,\enspace\text{and}\\
\forall\cQ&\in\MCL\tsp\ \forall\Q\in\cQ,\ j\in\{1,2\},\ 0<q<1:\quad \bigl(\cQ\setminus\{\Q\}\cup\bisect_{j,q}(\Q)\bigr)\in\MCL\tsp.
\end{align*}
\end{df}
\begin{df}[Skeleton]
Given a mesh $\cQ\in\MCL\tsp$ and $\Q=[x,x+\tilde x]\times[y,y+\tilde y]\in\cQ$, we denote the union of all vertical (resp. horizontal) element sides by 
\begin{align*}
\vsk(\Q) &\sei \{x,x+\tilde x\}\times[y,y+\tilde y],\\
\hsk(\Q) &\sei [x,x+\tilde x]\times\{y,y+\tilde y\},\\
\vsk(\cQ)&\sei \bigcup_{\Q\in\cQ}\vsk(\Q),\quad
\hsk(\cQ)\sei \bigcup_{\Q\in\cQ}\hsk(\Q).
\end{align*}
We call $\vsk$ the \emph{vertical skeleton} and $\hsk$ the \emph{horizontal skeleton}.
\end{df}
\begin{df}[Vertices and T-junctions]
For any mesh $\cQ\in\MCL\tsp$ and element $\Q=[x,x+\tilde x]\times[y,y+\tilde y]\in\cQ$, we define the set of \emph{vertices} 
\[\V(\Q)\sei\{x,x+\tilde x\}\times\{y,y+\tilde y\},\enspace\text{and}\enspace\V(\cQ)\sei\bigcup_{\Q\in\cQ}\V(\Q).\]
We denote as \emph{T-junction} each vertex that is in an element without being a vertex of it, 
\[\T(\Q)\sei \V(\cQ)\cap\Q\setminus\V(\Q),\enspace\text{and}\enspace\T(\cQ)\sei\bigcup_{\Q\in\cQ}\T(\Q).\]
Note that the above union is disjoint, i.e., for any T-junction $v\in\T(\cQ)$ there is a \emph{unique} element $\Q_v\in\cQ$ such that $v\in\T(\Q_v)$.
We distinguish horizontal and vertical T-junctions. A T-junction is called \emph{horizontal} if it is in a vertical side of the corresponding element, and \emph{vertical} if it is in a horizontal side,
\begin{align*}
\Th(\Q)&\sei\bigl\{v\in\T(\Q)\mid v\in \vsk(\Q)\bigl\},\\
\Tv(\Q)&\sei\bigl\{v\in\T(\Q)\mid v\in \hsk(\Q)\bigl\},\\
\Th(\cQ)&\sei\bigcup_{\Q\in\cQ}\Th(\Q),\quad
\Tv(\cQ)\sei\bigcup_{\Q\in\cQ}\Tv(\Q).
\end{align*}
Note that $\Th(\cQ)$ and $\Tv(\cQ)$ are disjoint and $\Th(\cQ)\cup\Tv(\cQ)=\T(\cQ)$.
\end{df}

\begin{df}
For any $v=(v_1,v_2)\in[0,M]\times[0,N]$, we define
\begin{alignat*}{2}
\mathsf X(v) &\sei \bigl\{&z\in\vsk&\mid z_2=v_2\bigr\}\cup\bigl\{ -\ceilfrac p2,\dots,-1,M+1,\dots, M+\ceilfrac p2\bigr\}\times\{y\},\\
\mathsf Y(v) &\sei \bigl\{&z\in\hsk&\mid z_1=v_1\bigr\}\cup\{x\}\times\bigl\{ -\ceilfrac q2,\dots,-1,N+1,\dots, N+\ceilfrac q2\bigr\}.
\end{alignat*}
\end{df}
\begin{rmrk}
The above-defined sets $\mathsf X(y)$ and $\mathsf Y(x)$ are sets of points, in contrast to the literature, where they are defined as sets of indices and referred to as ``global index sets'' \cite{mp2015} or ``global index vectors'' \cite{scott2012}.
\end{rmrk}

\begin{df}[T-junction extensions]
For any T-junction $v\in\T(\cQ)$, we define the T-junction extension as follows.
Consider $\mathsf X(v)$  to be ordered with respect to the first coordinate, then 
$\mathsf x_\mathsf{ext}(v)$  is defined as the unique set of $p_1+1$ consecutive elements of $\mathsf X(v)$  having the two elements of $\mathsf X(v)\cap\Q_v$ as the two middle entries.
We denote by $\conv(\mathsf x_\mathsf{ext}(v))$ the convex hull of these points.
Analogously, let $\mathsf Y(v)$ be ordered with respect to the second coordinate, and 
$\mathsf y_\mathsf{ext}(v)$  the unique set of $p_2+1$ elements of $\mathsf Y(v)$ having the two elements of $\mathsf Y(v)\cap\Q_v$ as the two middle entries, and $\conv(\mathsf y_\mathsf{ext}(v))$ the convex hull of these points.
The \emph{T-junction extension} of $v$ is defined as
\[\ext_\cQ(v)\sei\begin{cases}
\conv(\mathsf x_\mathsf{ext}(v))&\text{if }v\in\Th(\cQ),\\
\conv(\mathsf y_\mathsf{ext}(v))&\text{if }v\in\Tv(\cQ).
\end{cases}\]
\end{df}
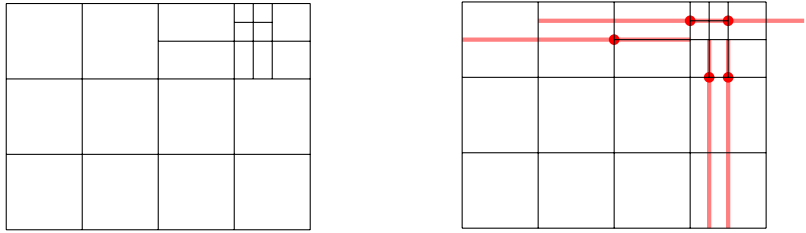
\begin{figure}[ht]
\centering
\begin{tikzpicture}
\draw (1,1) grid (5,4) (3,3.5)--(5,3.5) (4.25,3)--(4.25,4) (4.5,3)--(4.5,4) (4,3.75)--(4.5,3.75);
\end{tikzpicture}
\hspace{.1\textwidth}
\begin{tikzpicture}
\draw[ultra thick, red!50] (1,3.5)--(4,3.5)
                  (2,3.75)--(4.25,3.75)
                  (4.25,3.75)--(5.5,3.75)
                  (4.25,3.5)--(4.25,1)
                  (4.5,3.5)--(4.5,1);
\foreach \a/\b in {3/3.5, 4.25/3, 4.5/3, 4/3.75, 4.5/3.75}
  \fill[red] (\a,\b) circle (2pt);
\draw (1,1) grid (5,4) (3,3.5)--(5,3.5) (4.25,3)--(4.25,4) (4.5,3)--(4.5,4) (4,3.75)--(4.5,3.75);
\end{tikzpicture}
\caption{Example for T-junction extensions. The left figure show the considered mesh, and the right figure shows the same mesh with indicated T-junctions (red bullets) and the corresponding T-junction extensions (light red thick lines).}
\end{figure}
\begin{df}[B\'ezier mesh]\label{df: bezier mesh}%
Given a mesh $\cQ\in\MCL\tsp$, adding all T-junction extensions as actual edges to $\cQ$ yields the \emph{B\'ezier mesh}, also called \emph{extended T-mesh}. It represents the lowest-dimensional piecewise polynomial space that contains $\cB\tsp^\cQ$ and is used for mesh comparisons in this paper, see e.g.\ Figure~\ref{fig:SquareCorner_Mesh}.
\end{df}

\begin{df}[T-spline functions]
To each active node $v=(v_1,v_2)\in\V(\cQ)$, we associate  a local index vector $\mathsf x(v)\in\mathbb R^{p+2}$, which is obtained by taking the unique $p+2$ consecutive elements in $\mathsf X(v_2)$ having $v_1$ as their $\tfrac{p+3}2$-th (this is, the middle) entry. We analogously define $\mathsf y(v)\in\mathbb R^{q+2}$.

We associate to each active node $v\in\V(\cQ)$ a bivariate B-spline function defined as the product of the one-dimensional B-spline functions on the corresponding local index vectors, \[B_v(x,y,z)\coloneqq N_{\mathsf x(v)}(x) \cdot N_{\mathsf y(v)}(y).\]
Given a mesh $\cQ\in\MCL\tsp$, the associated set of T-spline functions is defined by
\[\cB\tsp^{\cQ}\sei\left\{B_v\mid v\in\V(\cQ)\right\}.\]
\end{df}
It is known from the literature that these functions are linearly independent if and only if there is no intersection between a horizontal and a vertical T-junction extension \cite{daveiga2012,sederberg2012}. Moreover, given a mesh $\cQ\in\MCL\tsp$ and a refinement $\tilde\cQ$ thereof, the corresponding spline spaces are only nested if each T-junction extension is either eliminated or unchanged \cite{scott2014}.
\begin{df}
For any T-junction $v\in\T(\cQ)$, we denote by \[\reftj(\cQ,v)\sei\cQ\setminus\{\Q_v\}\cup\bisect_{j,q}(\Q_v)\] the single-element refinement such that $v\notin\T(\reftj(\cQ,v))$, i.e., such that $v$ is not a T-junction anymore.
\end{df}
\begin{rmrk}
This refinement exists and is unique, and it is constructed as follows. The definition of T-junctions states that there is exactly one element $\Q_v\in\cQ$ such that $v\in\T(\Q_v)$. The location of $v$ on the boundary of $\Q_v$ uniquely defines bisection parameters $j$ and $q$ such that $v$ is a vertex of each children $\Q'\in\bisect_{j,q}(\Q_v)$.
\end{rmrk}
\begin{df}[Extension crossing, extension incompatibility]
For any mesh $\cQ\in\MCL\tsp$, we denote the set of \emph{extension-crossing T-junction pairs} by
\[\mathsf E(\cQ)\sei \left\{ (v,w)\in\Th(\cQ)\times\Tv(\cQ)\mid\ext_\cQ(v)\cap\ext_\cQ(w)\neq\emptyset\right\}.\]
For any mesh $\cQ\in\MCL\tsp$ and refinement $\tilde\cQ\in\MCL\tsp$, we define the set of \emph{extension-incompatible T-junctions} by
\[\mathsf C(\cQ,\tilde\cQ) \sei \left\{ v\in\T(\cQ)\cap\T(\tilde\cQ)\mid\ext_{\tilde\cQ}(v)\subsetneqq\ext_\cQ(v)\right\}.\]
\end{df}

\begin{alg}[\Little refinement for T-splines, \cite{scott2012}]\ 
\begin{algorithmic}
\setlength \baselineskip{14pt}
\REQUIRE {mesh $\cQ\in\MCL\tsp$, marked elements $\cM\subset\cQ$}
\STATE {$\tilde\cQ\sei\refthb(\cQ,\cM)$}
\REPEAT
\STATE {$\displaystyle 
\begin{aligned}
v_\text{refine}&\sei\operatorname*{argmin}_{v\in\T(\tilde\cQ)}\Bigl(\#\mathsf E(\reftj(\tilde\cQ,v)) + \#\mathsf C(\cQ,\reftj(\tilde\cQ,v))\Bigr)\\
\tilde\cQ&\sei\reftj(\tilde\cQ,v_\text{refine})\\[.75ex]
\end{aligned}$}
\UNTIL {$\mathsf E(\tilde\cQ)=\emptyset$} \AND {$\mathsf C(\cQ,\tilde\cQ)=\emptyset$}
\RETURN {$\reftsp\p(\cQ,\cM)\sei\bigl(\tilde\cQ,\cB\tsp^{\tilde\cQ}\bigr)$}
\end{algorithmic}
\end{alg}
\begin{rmrk}
The above algorithm does always finish, in the worst case yielding a tensor-product mesh.
\end{rmrk}

\begin{figure}[ht]
\captionsetup[subfigure]{labelformat=empty}
\centering
$\phantom\rightarrow$\quad
\begin{subfigure}{.25\textwidth}
\centering
\begin{tikzpicture}
\draw[white]      (4.125,3.625)--(4.125,4.5);
\fill[blue!50] (4,3.5) rectangle (4.25,3.75);
\draw (1,1) grid (5,4) (3,3.5)--(5,3.5) (4.25,3)--(4.25,4) (4.5,3)--(4.5,4) (4,3.75)--(4.5,3.75);
\end{tikzpicture}
\caption{$\cQ$\\$\phantom C$}
\end{subfigure}
\quad$\rightarrow$\quad
\begin{subfigure}{.25\textwidth}
\centering
\begin{tikzpicture}
\draw[white]      (4.125,3.625)--(4.125,4.5);
\draw (1,1) grid (5,4) (3,3.5)--(5,3.5) (4.25,3)--(4.25,4) (4.5,3)--(4.5,4) (4,3.75)--(4.5,3.75);
\draw[step=.125] (4,3.5) grid (4.25,3.75);
\end{tikzpicture}
\caption{$\tilde\cQ_0$\\$\phantom C$}
\end{subfigure}
\quad$\phantom\rightarrow$\quad
\begin{subfigure}{.25\textwidth}
\centering
\begin{tikzpicture}
\draw[very thick, red!50] (1,3.5)--(4,3.5)
                  (2,3.75)--(4.25,3.75)
                  (4.25,3.75)--(5.5,3.75)
                  (4.25,3.5)--(4.25,1)
                  (4.5,3.5)--(4.5,1)
                  (4.125,3.625)--(5,3.625)
                  (4.125,3.625)--(4.125,2)
                  (4.125,3.625)--(2,3.625)
                  (4.125,3.625)--(4.125,4.5);
\foreach \a/\b in {3/3.5, 4.25/3, 4.5/3, 4/3.75, 4.5/3.75, 4/3.625, 4.125/3.5, 4.25/3.625, 4.125/3.75}
  \fill[red] (\a,\b) circle (1.6pt);
\draw (1,1) grid (5,4) (3,3.5)--(5,3.5) (4.25,3)--(4.25,4) (4.5,3)--(4.5,4) (4,3.75)--(4.5,3.75);
\draw[step=.125] (4,3.5) grid (4.25,3.75);
\end{tikzpicture}
\caption{$\#\mathsf E(\tilde\cQ_0)=5$\\$\#\mathsf C(\tilde\cQ_0)=0$}
\end{subfigure}
\\[1em]$\rightarrow$\quad
\begin{subfigure}{.25\textwidth}
\centering
\begin{tikzpicture}
\draw[very thick, red!50] (1,3.5)--(4,3.5)
                  (2,3.75)--(4.25,3.75)
                  (4.25,3.75)--(5.5,3.75)
                  (4.25,3.5)--(4.25,1)
                  (4.5,3.5)--(4.5,1)
                  (4.125,3.625)--(5,3.625)
                  (4.125,3.625)--(4.125,2)
                  (4.125,3.625)--(2,3.625);
\foreach \a/\b in {3/3.5, 4.25/3, 4.5/3, 4/3.75, 4.5/3.75, 4/3.625, 4.125/3.5, 4.25/3.625}
  \fill[red] (\a,\b) circle (1.6pt);
\draw (1,1) grid (5,4) (3,3.5)--(5,3.5) (4.25,3)--(4.25,4) (4.5,3)--(4.5,4) (4,3.75)--(4.5,3.75);
\draw[step=.125] (4,3.5) grid (4.25,3.75) (4.125,3.75)--(4.125,4);
\end{tikzpicture}
\caption{$\#\mathsf E(\tilde\cQ_1)=2$\\$\#\mathsf C(\tilde\cQ_1)=0$}
\end{subfigure}
\quad$\rightarrow$\quad
\begin{subfigure}{.25\textwidth}
\centering
\begin{tikzpicture}
\draw[very thick, red!50] (1,3.5)--(4,3.5)
                  (2,3.75)--(4.25,3.75)
                  (4.25,3.75)--(5.5,3.75)
                  (4.25,3.5)--(4.25,1)
                  (4.5,3.5)--(4.5,1)
                  (4.125,3.625)--(5,3.625)
                  (4.125,3.5)--(4.125,1)
                  (4.125,3.625)--(2,3.625);
\foreach \a/\b in {3/3.5, 4.25/3, 4.5/3, 4/3.75, 4.5/3.75, 4/3.625, 4.125/3, 4.25/3.625}
  \fill[red] (\a,\b) circle (1.6pt);
\draw (1,1) grid (5,4) (3,3.5)--(5,3.5) (4.25,3)--(4.25,4) (4.5,3)--(4.5,4) (4,3.75)--(4.5,3.75);
\draw[step=.125] (4,3.5) grid (4.25,3.75) (4.125,3)--(4.125,4);
\end{tikzpicture}
\caption{$\#\mathsf E(\tilde\cQ_2)=0$\\$\#\mathsf C(\tilde\cQ_2)=0$}
\end{subfigure}
\quad$\phantom\rightarrow$\quad
\begin{subfigure}{.25\textwidth}
\centering
\begin{tikzpicture}
\draw (1,1) grid (5,4) (3,3.5)--(5,3.5) (4.25,3)--(4.25,4) (4.5,3)--(4.5,4) (4,3.75)--(4.5,3.75);
\draw[step=.125] (4,3.5) grid (4.25,3.75) (4.125,3)--(4.125,4);
\end{tikzpicture}
\caption{$\tilde\cQ_2$\\$\phantom C$}
\end{subfigure}
\caption{Example for the \little T-spline refinement. In the first step, the marked element is subdivided as for the THB-spline refinement. Second, the intersections of horizontal and vertical T-junction extensions are counted. Third, $\reftj$ is applied to a T-junction for which the number of extension crossings in the resulting mesh (elements of $\mathsf E(\tilde\cQ)$) is smallest, plus the term $\#\mathsf C(\cQ,\tilde\cQ)$ to ensure nesting of the resulting spline spaces. This third step is repeated until the sets $\mathsf E(\tilde\cQ)$ and $\mathsf C(\cQ,\tilde\cQ)$ are empty.}
\end{figure}
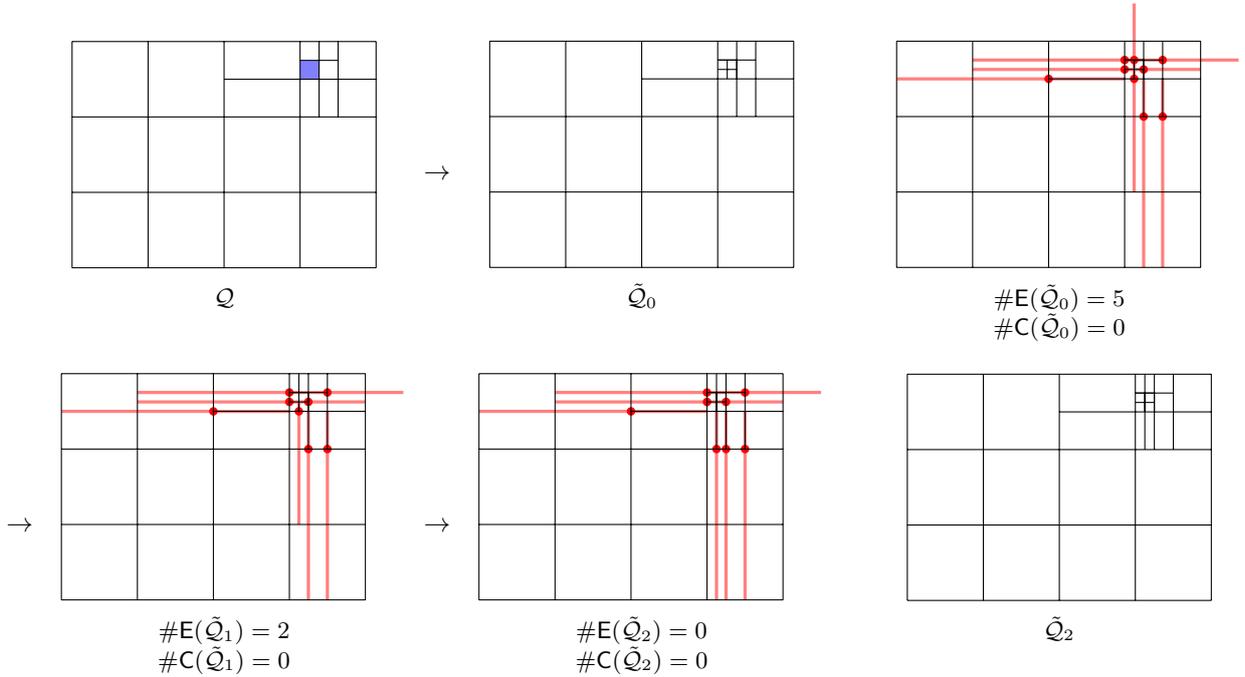

\subsubsection{\Much refinement for T-splines}
\begin{df}
For each level $k\in\nat$, we define the tensor-product mesh
\[\cQ_{k+1/2}\sei\Bigl\{[x-2^{-k-1},x]\times[y-2^{-k},y]\mid 2^{k+1} x\in\{1,\dots,2^{k+1} M\},2^k y\in\{1,\dots,2^k N\}\Bigr\}.\]
\end{df}
\begin{df}[intermediate children]
For $k\in\nat_0$ and $\Q\in\cQ_{k/2}$, we define
\[\subdiv^{1/2}(T)\sei \bigl\{T'\in\cQ_{(k+1)/2}\mid T'\subset T\bigr\}.\]
\end{df}

\begin{df}\label{df: reftsp}%
The T-spline refinement procedure reads
$\reftel(\cQ,\cM) \sei \bigl(\cQ\tsp,\cB\tsp\bigr)$, with
\begin{align*}
\tilde\cQ&\sei\cQ\setminus\cM\cup\subdiv^{1/2}(\cM)\\
\text{and }\quad\cB\tsp^{\tilde\cQ}&\sei\left\{B_v\mid v\in\V(\tilde\cQ)\right\}.
\end{align*}
\end{df}
\begin{df}[\Much refinement for T-splines]
We define for each $\Q\in\cQ$ the \emph{coarse neighbourhood}
\begin{align*}
\cN\tsppp(T)&\sei \left\{T'\in\cQ\cap\cQ_{\ell(\Q)-1/2}\mid \bigl(\operatorname{mid}(T)-\operatorname{mid}(T')\bigr)\le \mathbf D(\ell(\Q)) \right\},
\shortintertext{with} \mathbf D(k)&=\smash{ \begin{cases}
2^{-k/2}\left(\floorfrac p2+\tfrac12,\,\ceilfrac q2+\tfrac12\right)&\text{if $k$ is even,}\\
2^{-(k+1)/2}\left(\ceilfrac p2+\tfrac12,\,2\floorfrac q2+1\right)&\text{if $k$ is odd,}
\end{cases} }
\intertext{where $p$ and $q$ are the polynomial degrees of the B-splines in $x$- and $y$-direction, respectively.
Moreover, we define the \emph{closure}}
\clostsp\pp(\cM,\cQ)&\sei \smash{ \bigcup_{k=0}^{\max\ell(\cM)}\cN\tsppp^k(\cM), }
\intertext{and the extended refinement procedure}
\reftsp\pp(\cQ,\cM) &\sei \reftel\bigl(\cQ,\clostsp\pp(\cM,\cQ)\bigr).
\end{align*}
\end{df}

\begin{figure}[ht]
\centering
\includegraphics[width=.25\textwidth]{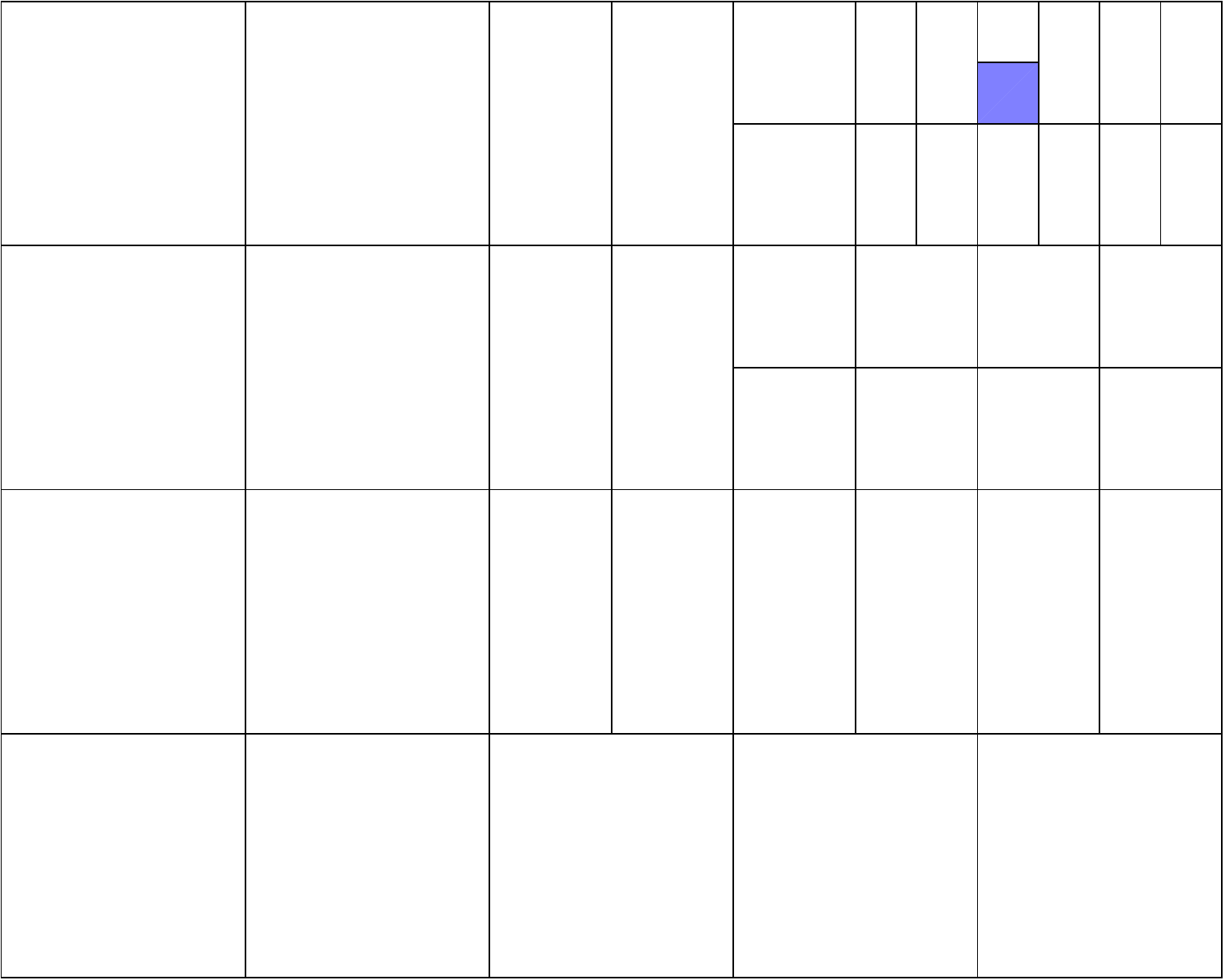}%
\quad\raisebox{.095\textwidth}{$\rightarrow$}\quad%
\includegraphics[width=.25\textwidth]{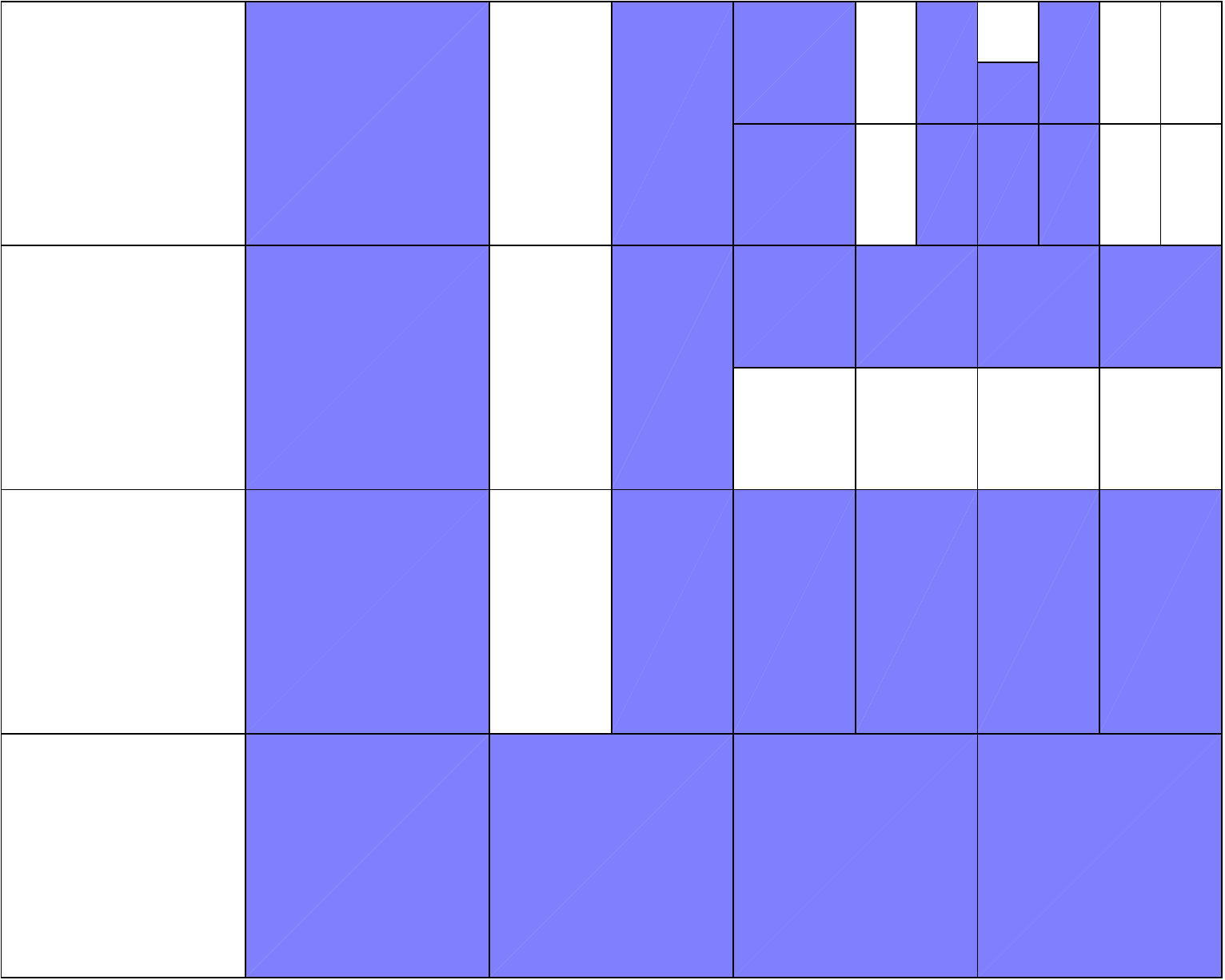}%
\quad\raisebox{.095\textwidth}{$\rightarrow$}\quad%
\includegraphics[width=.25\textwidth]{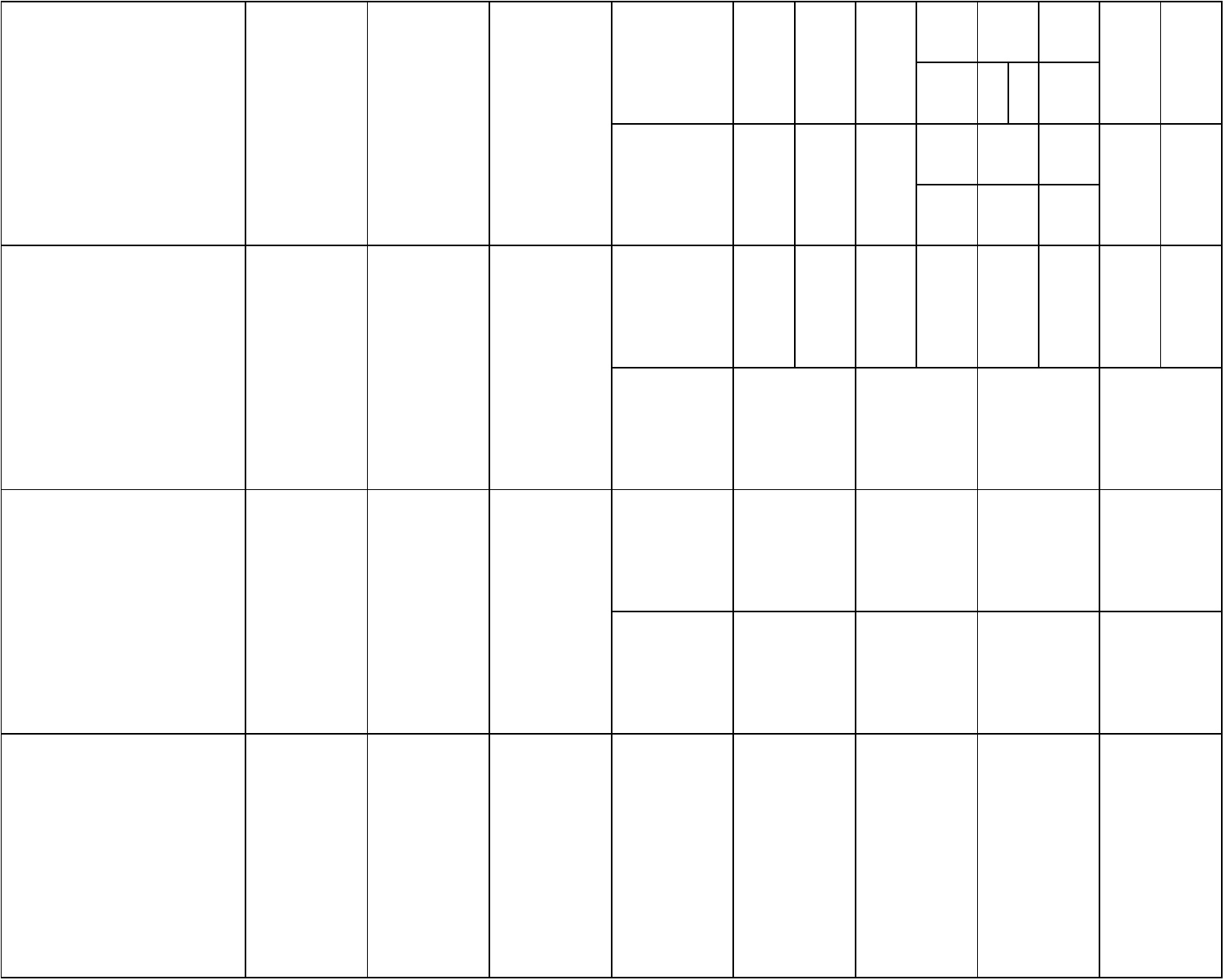}%
\\[1em]%
\includegraphics[width=.25\textwidth]{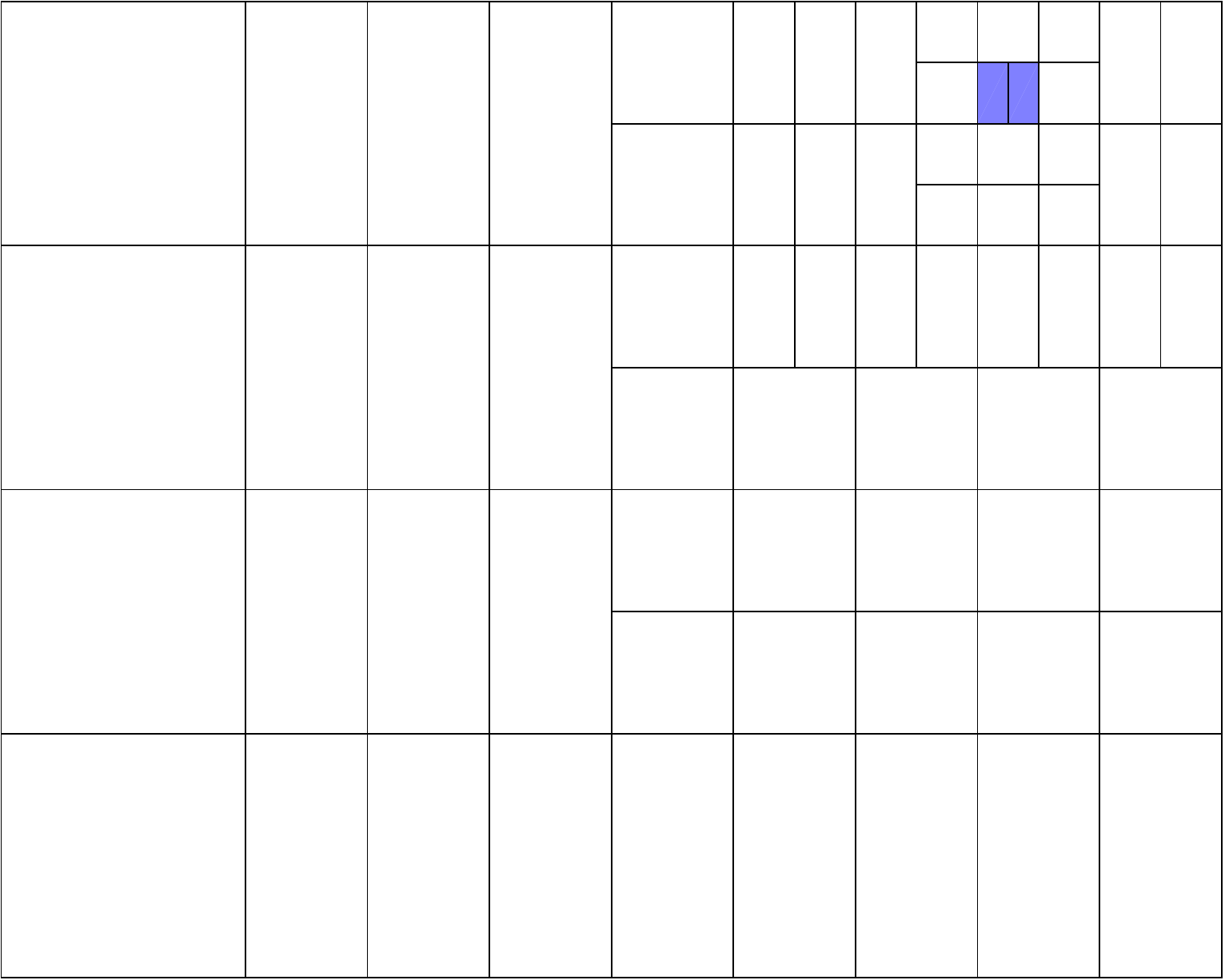}%
\quad\raisebox{.095\textwidth}{$\rightarrow$}\quad%
\includegraphics[width=.25\textwidth]{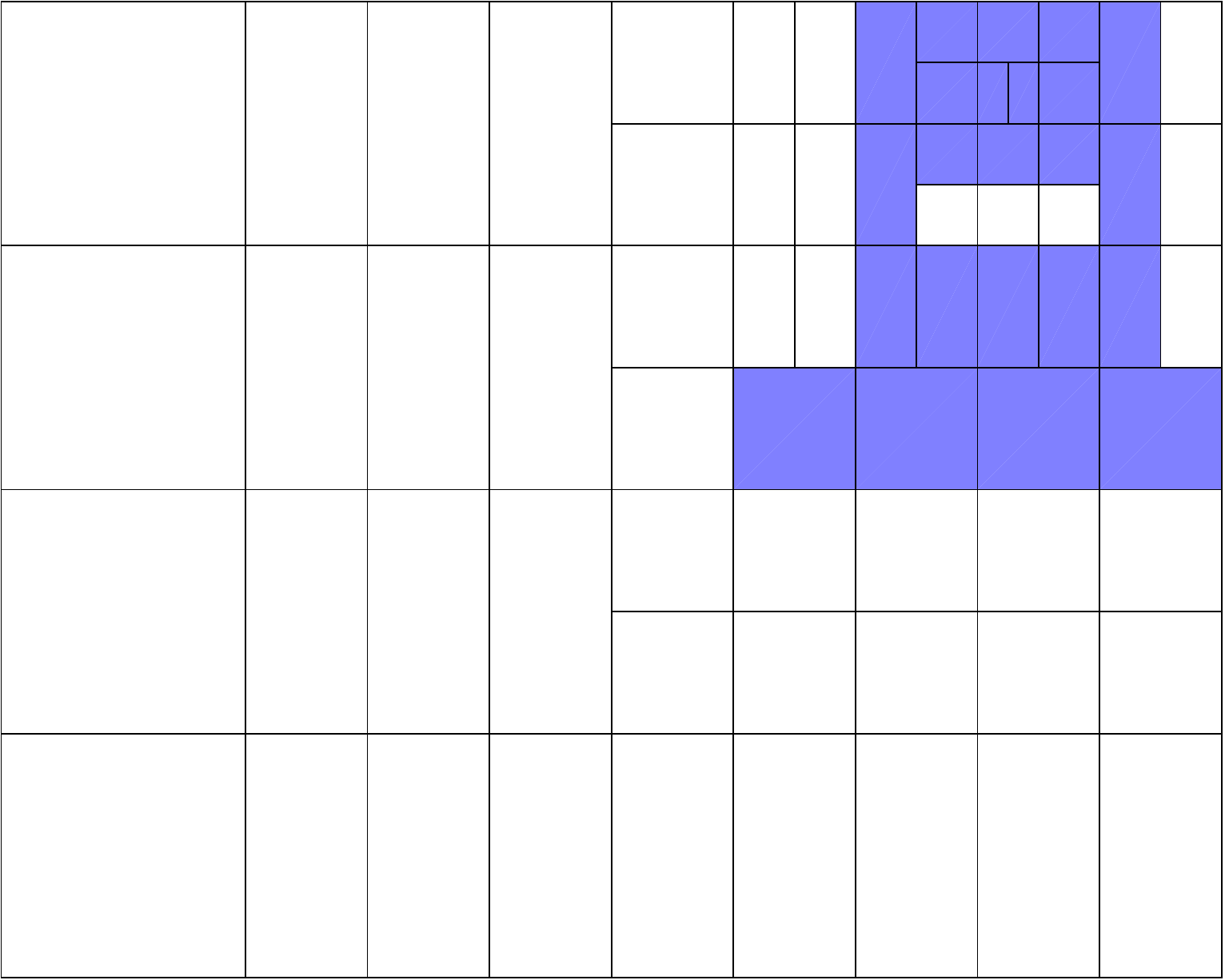}%
\quad\raisebox{.095\textwidth}{$\rightarrow$}\quad%
\includegraphics[width=.25\textwidth]{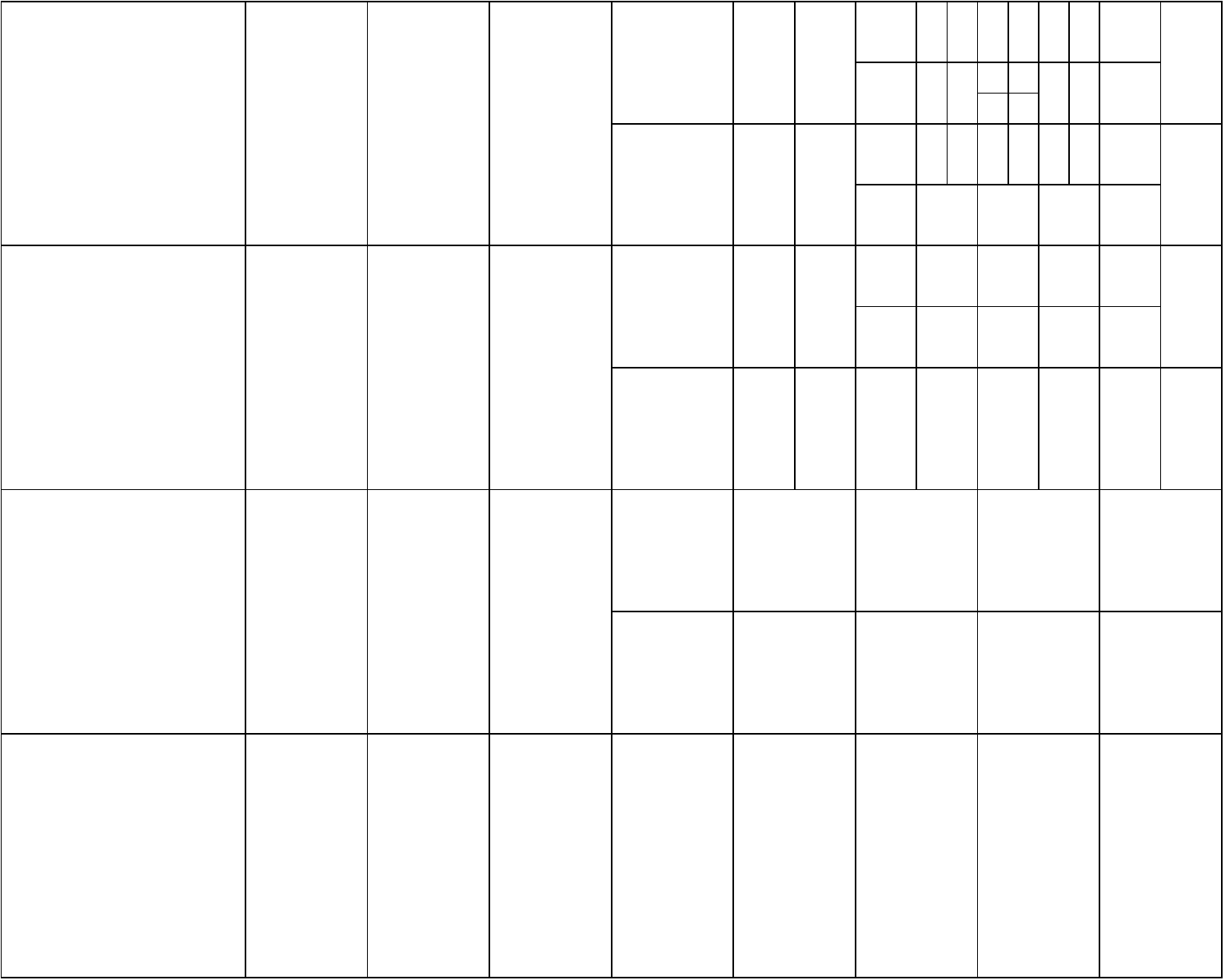}%
\caption{Example for the \much T-spline refinement. In order to subdivide the marked element as for the THB-spline refinement, the refinement routine $\reftsp\pp$ is applied twice.}
\end{figure}

\subsection{Theoretical Results}\label{sec: theo results}
As mentioned in the introduction, the use of T-splines and local mesh refinement faces fundamental difficulties. 
Given a mesh $\cQ\in\MCL\tsp$ and a refinement $\cQ'$ of $\cQ$, it is not clear in general that 
the T-spline functions that should serve as a spline basis
are in fact linearly independent, and, even if they are, that the new spline space $\Span\cB\tsp^{\cQ'}$ is a superspace of the preceeding space $\Span\cB\tsp^\cQ$.
Both refinement procedures $\reftsp\p$ and $\reftsp\pp$ overcome these problems.

Each of the refinement procedures $\refine\in\bigl\{\refthb\p,\refthb\pp,\reftsp\pp\bigr\}$ satisfies the following properties \cite{bgmp2016,mp2015}:
\begin{itemize}
\item \emph{Linear complexity.} There exists a constant $C$ that depends only on the polynomial degree of the B-splines in use, such that any sequence of meshes $\cQ_0,\cQ_1,\dots,\cQ_J$ with \[\cQ_j=\refine(\cQ_{j-1},\cM_{j-1}),\quad\cM_{j-1}\subseteq\cQ_{j-1}\quad\text{for}\enspace j\in\{1,\dots,J\}\] satisfies
\[\left|\cQ_J\setminus\cQ_0\right|\ \le\ C\sum_{j=0}^{J-1}|\cM_j|\ .\]
This leads to linear complexity in terms of degrees of freedom. However, it does in general not reflect the complexity with respect to computing time or memory.
\item \emph{Bounded overlay.} For $\cQ_a$ and $\cQ_b$ being two meshes generated by $\cQ_0$ and the successive use of $\refine$, there exists a common refinement of $\hat\cQ$ of $\cQ_a$ and $\cQ_b$ such that
\[\#\hat\cQ+\#\cQ_0\le\#\cQ_a+\#\cQ_b.\]
\end{itemize}
Moreover, each of these refinement procedures has a natural generalization to higher dimensions \cite{giannelli2012,morgenstern2015}, and due to the hierarchical construction of the mesh classes, they provably provide a stable shape regularity in the sense of bounded aspect ratios of the mesh elements and local quasi-uniformity of the mesh \cite[Lemma~2.14]{mp2015}.
Such analysis or higher-dimensional version is currently unavailable for $\reftsp\p$.

\section{Model Problems and Discretization}\label{sec: model problem}
This section describes the two model problems that are used for our tests. We will formulate both problems in the weak (variational) form and skip their derivation from the original PDEs. The latter are, for the Poisson problem, seeking $u\in C^2(\overline\Omega)$ such that \[-\Delta u=f\text{ in }\Omega,\quad \tfrac{\partial u}{\partial\neunormal}=\neuBC\text{ on }\neubound\quad\text{and}\quad u|_\dirbound=\dirBC\text{ on }\dirbound,\]
and for the problem of linear elasticity, seeking $u\in C^2(\overline\Omega)$ such that
\[-\operatorname{div}\sigma(u)=f\text{ in }\Omega,\quad \langle\neunormal,\sigma(u)\rangle=\neuBC\text{ on }\neubound\quad\text{and}\quad u|_\dirbound=\dirBC\text{ on }\dirbound,\]
using the notation explained below.
\subsection{Poisson problem} \label{sec: Poisson problem}
\paragraph{Data}
Let $\Omega\subset\reell^d$ be an open, connected and bounded Lipschitz domain. Let the Dirichlet boundary $\dirbound\subset\partial\Omega$ be closed and let each connectivity component of $\dirbound$ be of positive measure. Set the Neumann boundary $\neubound\sei\partial\Omega\setminus\dirbound$ and the corresponding outer normal vector $\neunormal:\neubound\to\reell^d$. Let $\dirBC\in L^2(\dirbound)$ and 
\[H^1_0(\Omega)\sei\left\{w\in H^1(\Omega)\mid w|_\dirbound=0\text{ a.e.\ in }\dirbound\right\}.\]
\paragraph{Problem}
Find $u\in H^1(\Omega)$ such that 
\begin{equation}\label{eq: poisson weak}
\begin{aligned}
\int_\Omega\langle\nabla u,\nabla v\rangle \d x &= \int_\Omega \RHS v\d x+\int_\neubound \neuBC v\d s\quad\text{for all }v\in H^1_0,\\
u|_\dirbound&=\dirBC\quad\text{a.e.\ on }\dirbound.
\end{aligned}
\end{equation}
\paragraph{Discretization} 
Given a basis $\cB$ of a finite-dimensional function space $\hat\cB\sei\Span\cB$ and
\[\hat\cB_0\sei\left\{w\in\hat\cB\mid w|_\dirbound=0\text{ a.e.\ in }\dirbound\right\},\]
we seek the Galerkin solution $\hat u\in\hat\cB$ satisfying
\begin{equation}\label{eq: poisson weak discrete}
\begin{aligned}
\int_\Omega\langle\nabla \hat u,\nabla v\rangle \d x &= \int_\Omega \RHS v\d x+\int_\neubound \neuBC v\d s\quad\text{for all }v\in \hat\cB_0,\\
\hat u|_\dirbound&=I_{\hat\cB}(\dirBC)\quad\text{on }\dirbound,
\end{aligned}
\end{equation}
where $I_{\hat\cB}(\dirBC)\in\hat\cB$ is an interpolation of $\dirBC$.
We set $\hat u_0\sei \hat u -I_{\hat\cB}(\dirBC)\in\hat\cB_0$ and reformulate the above problem to finding $\hat u_0\in\hat\cB_0$ such that 
\begin{equation}\label{eq: poisson weak discrete homBC}
\int_\Omega\langle\nabla \hat u_0,\nabla v\rangle \d x = \int_\Omega \RHS v\d x+\int_\neubound \neuBC v\d s
- \int_\Omega\langle\nabla I_{\hat\cB}(\dirBC),\nabla v\rangle \d x\quad\text{for all }v\in \hat\cB_0.
\end{equation}
Since both left and right side of \eqref{eq: poisson weak discrete} are linear in $v$, it suffices to have the above equation fulfilled for all basis functions $v\in\cB_0=\{v_1,\dots,v_n\}=\cB\cap\hat\cB_0$ that are zero on the boundary. Since the right-hand side is also linear in $\hat u_0$, and $\hat u_0\in\hat\cB_0$ is a linear combination of these basis functions, \eqref{eq: poisson weak discrete homBC} is equivalent to finding a vector $U=(u_1,\dots,u_n)\in\reell^n$ such that
\begin{equation}\label{eq: poisson discrete}
\underbrace{\Bigl(\int_\Omega\langle\nabla v_i,\nabla v_j\rangle \d x\Bigr)_{1\le i,j\le n}}_{A\in\reell^{n\times n}}\ \cdot\  U = \underbrace{\Bigl(\int_\Omega \RHS v_i\d x+\int_\neubound \neuBC v_i\d s
- \int_\Omega\langle\nabla I_{\hat\cB}(\dirBC),\nabla v_i\rangle \d x
\Bigr)_{1\le i\le n}}_{B\in\reell^n},
\end{equation}
with $\hat u=\sum_{i=1}^nu_iv_i + I_{\hat\cB}(\dirBC)$.
We call $A$ the \emph{stiffness matrix} and $B$ the \emph{load vector}.
\paragraph{Error estimator} The Adaptive Algorithm (explained below in Section~\ref{sec: adaptive loop}) is controlled by a standard residual local error estimator $\eta:\cQ\to\reell$ (see e.g.\ \cite{buffa2015a} for an application with THB-splines). Given the Galerkin solution $\hat u\in\hat\cB$, it is defined by
\[\eta_\cQ(\Q)\sei \Bigl(
h_\Q^2\left\|\Delta \hat u+\RHS\right\|_\Q^2
+
\sum_{E\in\mathcal E(\Q)}h_E\left\|R_E(\hat u)\right\|_E^2
\Bigr)^{1/2},\]
where $\mathcal E(\Q)$ is the set of edges of $\Q$, $h_\Q$ the diameter of $\Q$, and $h_E$ the length (the 1D Lebesgue measure) of the edge $E$. The notation $\left\|\bullet\right\|_A$ abbreviates the $L^2$-norm $\left\|\bullet\right\|_{L^2(A)}$.
The \emph{edge residual} $R_E(\hat u)$ is defined by 
\[R_E(\hat u)\sei\begin{cases}
\tfrac12\,\big[\negthickspace\big[\tfrac{\partial \hat u}{\partial\nu_E}\big]\negthickspace\big]_E&\text{if $E$ is an interior edge,}\\[.5ex]
\neuBC-\tfrac{\partial \hat u}{\partial\nu_E}&\text{if $E$ is a boundary edge.}
\end{cases}
\]
For any interior edge $E=\Q\cap\Q'$, the notation $[\![\bullet]\!]_E\sei\bullet|_\Q-\bullet|_{\Q'}$ describes the jump along the edge $E$. Note that in all four methods  this paper accounts for, none of the spline basis functions have jumps in their derivatives, and the same holds for the discrete solution $\hat u$. Provided that the Neumann boundary condition is met exactly (e.g.\ in the case $\neuBC=0$), the above error estimator hence reduces to
\[\eta_\cQ(\Q)\sei 
h_\Q \left\|\Delta \hat u+\RHS\right\|_\Q
.\]

\subsection{Linear elasticity}
\paragraph{Data}
Let $\Omega$, $\dirbound$, $\neubound$, $\neunormal$, $\dirBC$ as above. For $u\in H^1(\Omega,\reell^d)$, we define
\[\varepsilon(u)\sei\bigl(\tfrac12(\partial_iu_j+\partial_ju_i)\bigr)_{1\le i,j\le d}
\quad\text{and}\quad
\sigma(u)_{ij} \sei \sum_{1\le k,\ell\le d}C_{ijk\ell}\,\varepsilon(u)_{k\ell}\enspace\text{for }1\le i,j\le d.\]
We call $u$ the \emph{displacement}, $\varepsilon(u)$ the \emph{strain tensor} and $\sigma(u)$ the \emph{stress tensor}, and $C$ is some positive definite fourth order tensor that describes material properties. 
\paragraph{Problem}
Find $u\in H^1(\Omega)$ such that 
\begin{equation}\label{eq: elasticity weak}
\begin{aligned}
\int_\Omega\langle\sigma(u),\varepsilon(v)\rangle \d x &= \int_\Omega \RHS v\d x+\int_\neubound \neuBC v\d s\quad\text{for all }v\in H^1_0(\reell^d),\\
u|_\dirbound&=\dirBC\quad\text{a.e.\ on }\dirbound.
\end{aligned}
\end{equation}
\paragraph{Discretization} 
Given a basis $\cB$ of a finite-dimensional function space $\hat\cB\sei\Span\cB$ and $\hat\cB_0$ as above, 
we seek the Galerkin solution $\hat u\in\hat\cB$ satisfying
\begin{equation}\label{eq: elasticity weak discrete}
\begin{aligned}
\int_\Omega\langle\sigma(\hat u),\varepsilon(v)\rangle \d x &= \int_\Omega \RHS v\d x+\int_\neubound \neuBC v\d s\quad\text{for all }v\in \hat\cB_0,\\
\hat u|_\dirbound&=I_{\hat\cB}(\dirBC)\quad\text{on }\dirbound.
\end{aligned}
\end{equation}
Again, we set $\hat u_0\sei \hat u -I_{\hat\cB}(\dirBC)\in\hat\cB_0$ and reformulate the above problem to finding $\hat u_0\in\hat\cB_0$ such that 
\begin{equation}\label{eq: elasticity weak discrete homBC}
\int_\Omega\langle\sigma(\hat u_0),\varepsilon(v)\rangle \d x = \int_\Omega \RHS v\d x+\int_\neubound \neuBC v\d s
- \int_\Omega\langle\sigma(I_{\hat\cB}(\dirBC)),\varepsilon(v)\rangle \d x\quad\text{for all }v\in \hat\cB_0.
\end{equation}
Analogously to the derivation for the Poisson problem above, we compute the Galerkin solution by solving the equation
\begin{equation}\label{eq: elasticity discrete}
\underbrace{\Bigl(\int_\Omega\langle\sigma(v_i),\varepsilon(v_j)\rangle \d x\Bigr)_{1\le i,j\le n}}_{A\in\reell^{n\times n}}\ \cdot\  U = \underbrace{\Bigl(\int_\Omega \RHS v_i\d x+\int_\neubound \neuBC v_i\d s - \int_\Omega\langle\sigma(I_{\hat\cB}(\dirBC)),\varepsilon(v_i)\rangle \d x\Bigr)_{1\le i\le n}}_{B\in\reell^n}.
\end{equation}
and setting $\hat u=\sum_{i=1}^nu_iv_i + I_{\hat\cB}(\dirBC)$.
\paragraph{Error estimator} Given the Galerkin solution $\hat u\in\hat\cB$, we use the local error estimator described in \cite{verfuerth1999}, which is defined by
\[\eta_\cQ(\Q)\sei \Bigl(
h_\Q^2\left\|\operatorname{div}\sigma(\hat u)+\RHS\right\|_\Q^2
+
\sum_{E\in\mathcal E(\Q)}h_E\left\|R_E(\hat u)\right\|_E^2
\Bigr)^{1/2},\]
where the \emph{edge residual} $R_E(\hat u)$ is defined by 
\[R_E(\hat u)\sei\begin{cases}
\tfrac12\left[\![\langle\nu_E,\sigma(\hat u)\rangle]\!\right]_E&\text{if $E$ is an interior edge,}\\
\neuBC-\langle\nu_E,\sigma(\hat u)\rangle&\text{if $E$ is a boundary edge.}
\end{cases}
\]

\section{Adaptive Algorithm}\label{sec: adaptive loop}
\subsection{Adaptive Loop}
The Adaptive Algorithm is an iterative procedure that consist of the steps
\[\tikz[remember picture] \coordinate (1);\textsf{SOLVE}\to\textsf{ESTIMATE}\to\textsf{MARK}\to\textsf{REFINE}\tikz[remember picture] \coordinate (2);\]
which are described as follows.
\begin{tikzpicture}[overlay, remember picture]
\foreach \a in {.075} {%
\draw (2) ++(2*\a,\a)--++(\a,0) arc (90:-90:2*\a) coordinate (3) (1) ++(-2*\a,\a) coordinate (4); 
\draw[decoration={markings,mark=at position 0.01 with {\arrow[scale=1.5]{<}}},
    postaction={decorate}] (4)--++(-\a,0) arc (90:270:2*\a) -- (3); }
\end{tikzpicture}
\begin{description}
\item[\sf SOLVE:] Given a finite-dimensional function space, compute a Galerkin approximation of the solution of the PDE.
\item[\sf ESTIMATE:] Compute local estimates for the error, i.e., the difference of approximate and exact solution.
\item[\sf MARK:] Based on these local estimates, select mesh elements $\cM\subseteq\cQ$ for refinement.
\item[\sf REFINE:] Refine the mesh $\cQ$ and construct the new discrete function space $\cB$.
\end{description}
Due to their dependence on the particular problem, the modules \textsf{SOLVE} and \textsf{ESTIMATE} have been defined above, for the two problems considered. For the module \textsf{REFINE}, we consider four variants, which have been outlined in Section~\ref{sec: mesh refinement}.

\subsection{Marking Strategies}\label{subsec: marking}
Given the estimated local errors $\{\eta_\Q\mid \Q\in\cQ\}\subset\reell$ and a marking parameter $\theta\in[0,1]$, which is chosen manually, the following strategies are commonly used for the step \textsf{MARK}.
\begin{itemize}
\item \emph{Quantile marking:} Let $\cQ=\{\Q_1,\dots,\Q_K\}$  and $\eta_\cQ(\Q_1)\ge\dots\ge\eta_\cQ(\Q_K)$, then $\cM=\{\Q_1,\dots,\Q_k\}$ with $k\approx\theta K$.
\item \emph{D\"orfler marking:} Let $\cQ=\{\Q_1,\dots,\Q_K\}$  and $\eta_\cQ(\Q_1)\ge\dots\ge\eta_\cQ(\Q_K)$, then $\cM=\{\Q_1,\dots,\Q_k\}$ with 
\[\sum_{j=1}^{k-1}\eta_\cQ(\Q_j)<\theta\sum_{\Q\in\cQ}\eta_\cQ(\Q)\quad\text{and}\quad\sum_{j=1}^k\eta_\cQ(\Q_j)\ge\theta\sum_{\Q\in\cQ}\eta_\cQ(\Q).\]
\item \emph{Maximum marking:} $\cM\sei\left\{\Q\in\cQ\mid \eta_\cQ(\Q)\ge\theta\max_{\tilde \Q\in\cQ}\eta_\cQ(\tilde \Q)\right\}$.
\end{itemize}

\subsection{Optimality of the Adaptive Algorithm} 
In the case of Quantile marking, the authors are not aware of theoretical results that ensure optimality of the convergence rates.
If the Adaptive FEM is applied with D\"orfler marking, the sequence of discrete solutions $u_1\in\hat\cB_1\in\BCL$, $u_2\in\hat\cB_2\in\BCL$, \dots has the best convergence rate (w.r.t.\ degrees of freedom) that is possible in the class $\BCL$ of discrete function spaces \cite{binev2004,stevenson2007,cascon2008,carstensen2014}.
For a modified version of Maximum marking, the discrete solutions are proven to be instance-optimal \cite{diening2015}. This means that in each step, the error of the discrete solution $u\in\hat\cB\in\BCL$ is bounded by a constant times the smallest error of the discrete solutions of all function spaces $\hat\cB'\in\BCL$ with a comparable (or smaller) number of degrees of freedom. This is an even stronger result than the rate optimality described above, however \cite{diening2015} accounts only for the Poisson problem with homogeneous Dirichlet boundary condition, i.e., $\neubound=\emptyset$ and $\dirBC=0$, and the authors are only aware of a generalization for the nonconforming Crouzeix-Raviart AFEM and the Stokes equation \cite{kreuzer2015}.

\section{Numerical Experiments}\label{sec: comparison}
In this section, the mesh refinement strategies from Section~\ref{sec: mesh refinement} are compared numerically. Hence, T-splines are compared with THB-splines and \little with \much refinement. In addition to achievable convergence rates and the mesh grading, the comparison includes the numerical properties of the stiffness matrix as its sparsity and condition number. 
To clearly point out differences between the refinement strategies, the first example is designed as a worst case scenario and does not correspond to a physical problem. The second and third example are well-established benchmark problems in the context Adaptive Finite Element Methods \cite{doerfel2010,Vuong2011,Hennig2016}, including the Poisson problem and linear elasticity with given analytical solutions. In all examples and for all refinement strategies, cubic B-spline basis functions are used.
\begin{figure}
 \centering
 \includegraphics[scale=0.80]{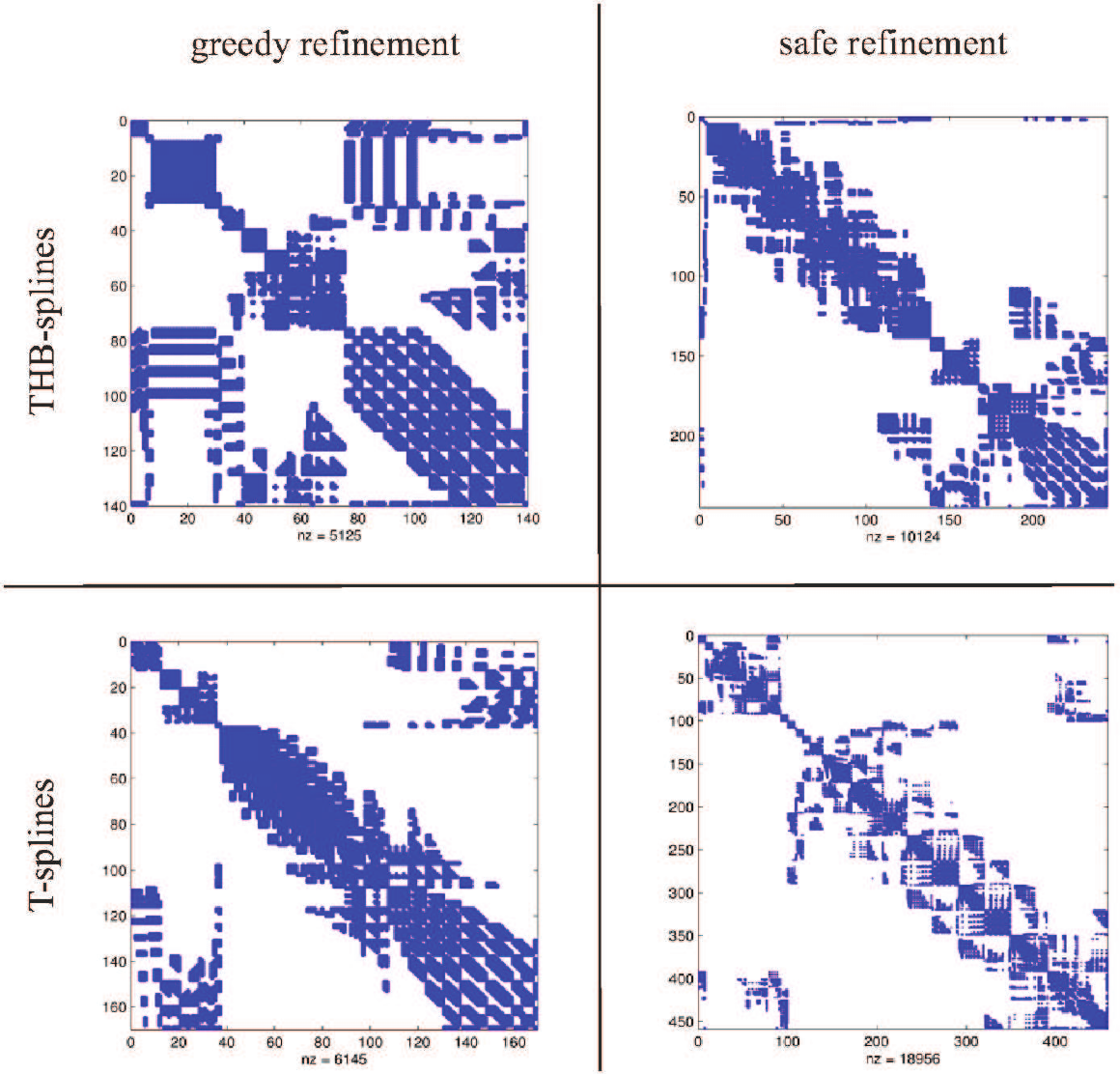}
 \caption{Worst case scenario: The sparsity patterns of the stiffness matrices after six refinement steps are illustrated. Especially the \little THB-spline refinement results in a dense stiffness matrix.}
 \label{fig:SquareCorner_Spar}
\end{figure}

\subsection{Worst case scenario}
\begin{figure}
 \centering
 \includegraphics[scale=0.80]{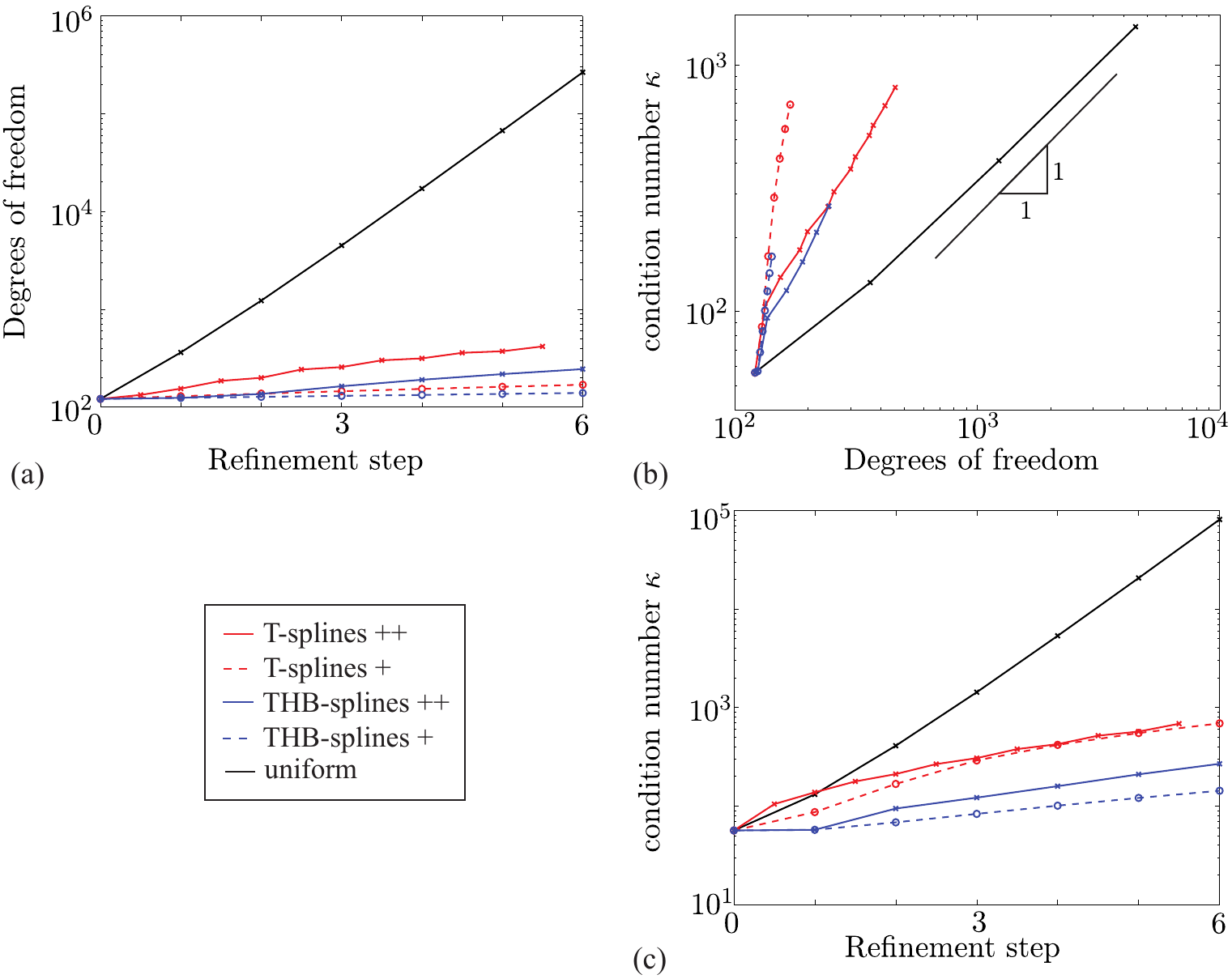}
 \caption{Worst case scenario: The relations between the total number of degrees of freedom, the condition number of the stiffness matrix and the refinement steps are illustrated.}
 \label{fig:SquareCorner_Conv}
\end{figure}
In this example, an initial square mesh with 64 elements is locally refined in the lower left corner, where only one element is marked for refinement in each refinement loop. The resulting B\' ezier meshes are presented in Figure~\ref{fig:SquareCorner_Mesh}. It can be seen that the \little THB-spline refinement does only refine the marked element whereas the \much refinement routines extend the refinement region. Also the \little T-spline refinement has to insert additional control points to ensure analysis-suitability. The total number of degrees of freedom (DOF) is plotted against the refinement steps in Figure~\ref{fig:SquareCorner_Conv}~(a) to illustrate this behaviour. 

The locality of the refinement comes at the cost of an increased interaction between differently scaled basis functions (cf. Section \ref{sec:Truncated Hierarchical B-splines}) in the case of \little THB-spline refinement. In this example, basis functions from the coarsest level interact with basis functions of the finest level. 
This leads to the occurrence of quasi-dense rows and columns and the loss of any band structure in the stiffnes matrix,
as it can be seen in Figure~\ref{fig:SquareCorner_Spar}. The other refinement routines do not produce anomalies in their sparsity patterns. 

The local mesh refinement also influences the behaviour of the condition number of the stiffness matrix. 
Gahalaut et al. \cite{Gahalaut2014} analyzed these condition numbers for NURBS-based isogeometric discretizations, showing that the condition number increases linearly with respect to degrees of freedom. This is also reflected in all our experiments. As expected, we observe for all kinds of local refinement that the condition numbers grow at higher rates, see Figure~\ref{fig:SquareCorner_Conv}~(b). The rate is apparently independent of the type (T- or THB-splines) but does depend on the locality of refinement (\little or \much), and thus on the grading of the mesh.
However, if the condition numbers are compared with respect to the refinement step (cf. Figure~\ref{fig:SquareCorner_Conv}~(c)), the \much THB-spline refinement produces higher condition numbers than the \little one, and the T-splines higher condition numbers than the THB-splines. 
This shows that the number of additional DOF per refinement step can has a dominant influence on the condition number. Hence, for a clear comparison, the condition number has to be compared with respect to a quantity of main interest. For this reason the numerical error of the solution will be plotted over the condition number in the following examples.
We emphasize that this discussion disregards appropriate preconditioning, but is beyond the scope of this paper.
\begin{figure}
 \centering
 \includegraphics[scale=0.80]{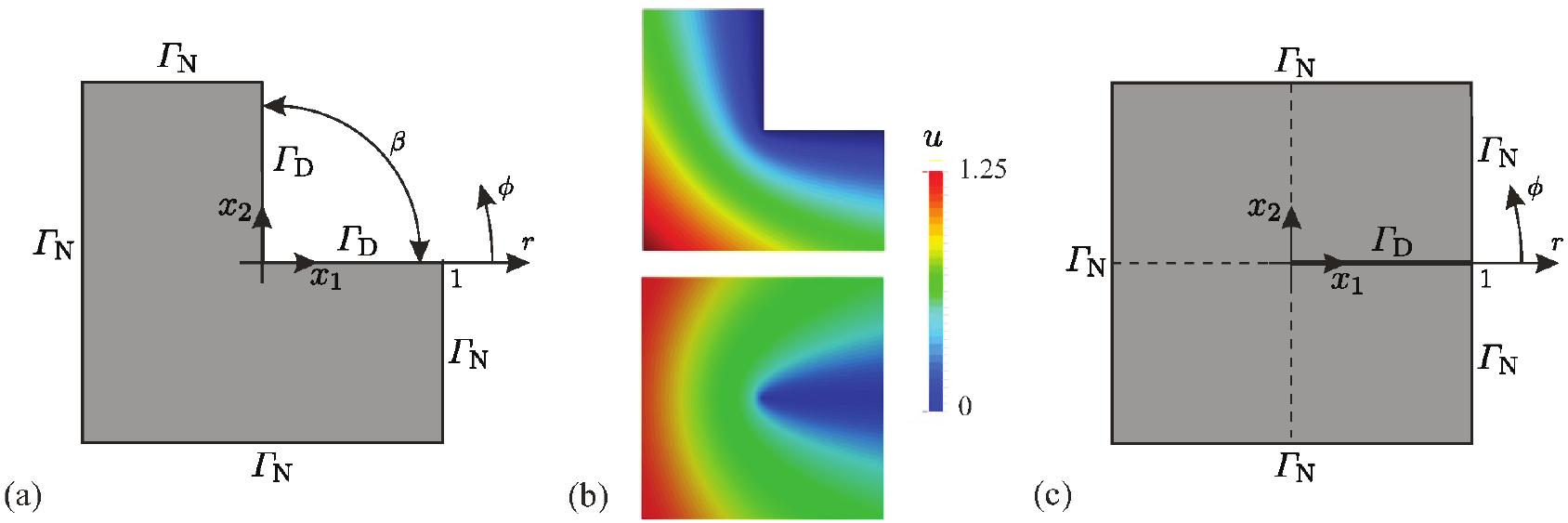}
 \caption{Poission problem: Domain and boundary conditions for \textbf{(a)} the L-shape and \textbf{(c)} the slit domain as well as \textbf{(b)} the corresponding analytical solutions.}
 \label{fig:HeatConductionAufbau_1}
\end{figure}
\begin{figure}
 \centering
 \includegraphics[scale=0.80]{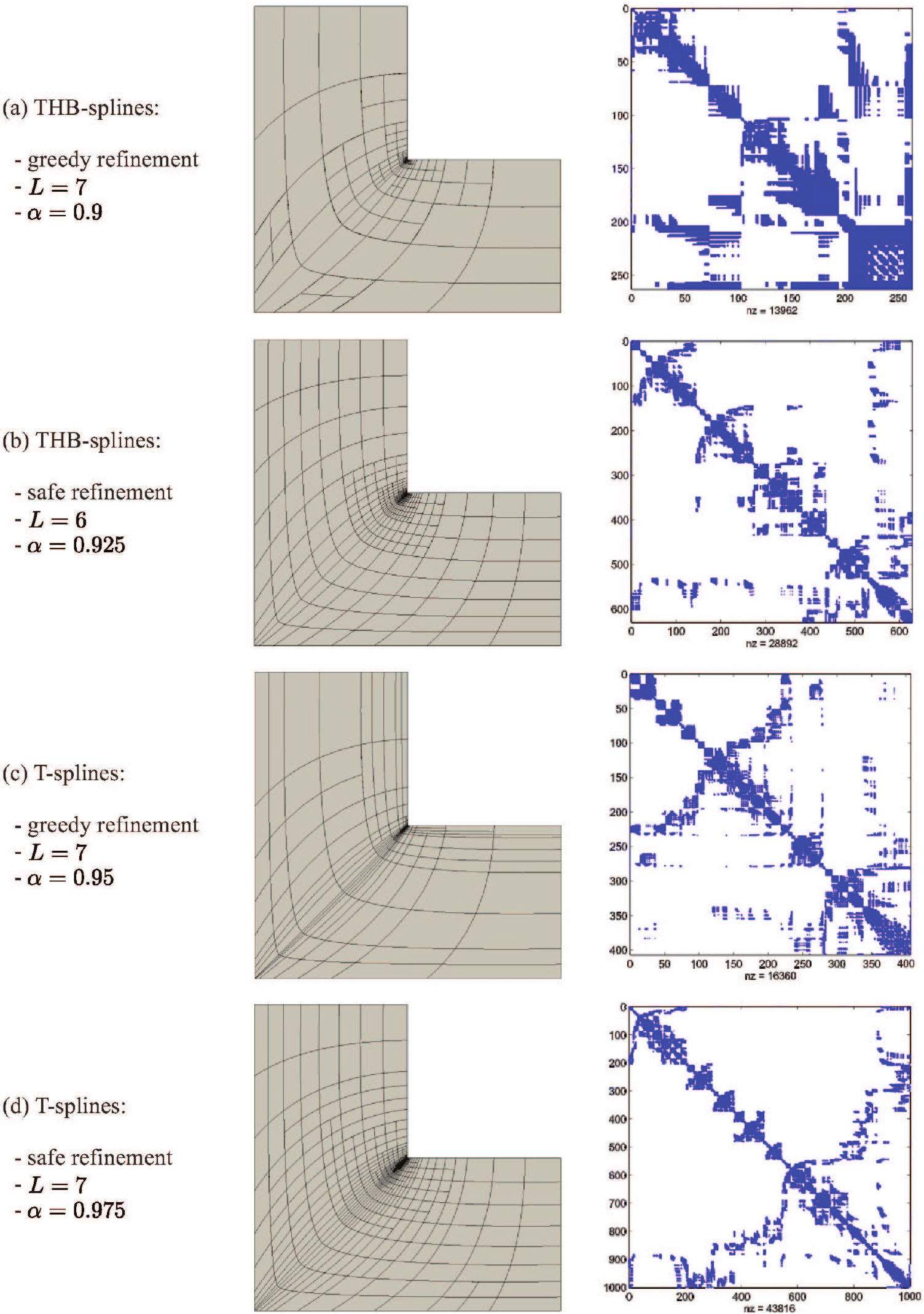}
 \caption{L-shape: The marking parameters $\alpha$, the B\' ezier meshes and the sparsity patterns of the stiffness matrices after $L$ refinement steps for all \textbf{(a)}-\textbf{(d)} refinement strategies. The \much refinement strategies result in well graded meshes, the \little refinement strategies in more unstructured meshes. Again, the \little THB-spline refinement creates the stiffness matrix with the highest density and interaction.}
 \label{fig:HeatConductionLShape_1}
\end{figure}
\begin{figure}
 \centering
 \includegraphics[scale=0.7]{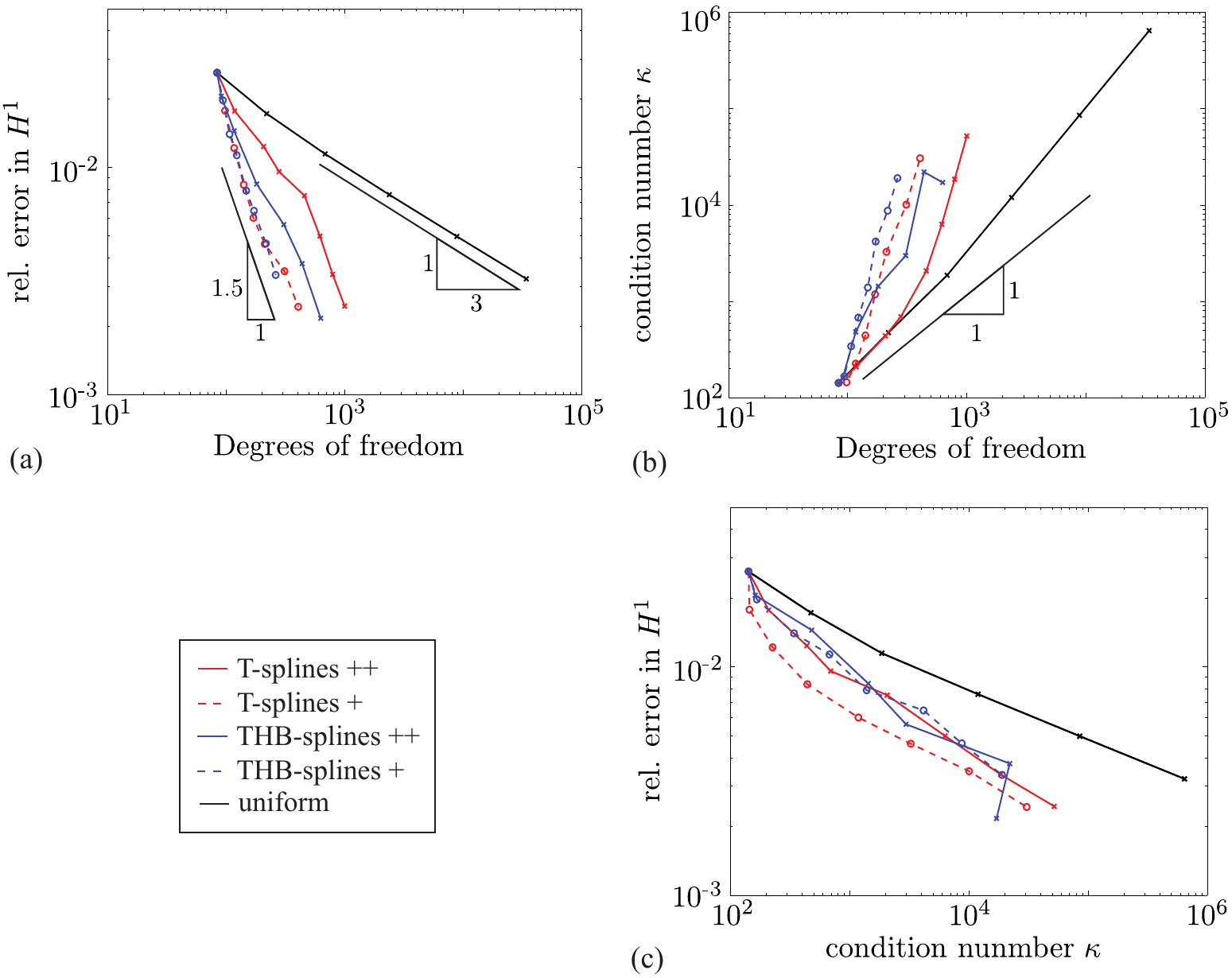}
 \caption{L-shape: The convergence rates as well as the relations between the condition number of the stiffness matrix, the numerical error of the solution and the total number of degrees of freedom are illustrated. \textbf{(a)} - All refinement strategies converge with the expected convergence rate $k=1.5$ in the asymptotic range.}
 \label{fig:HeatConductionLShape_2}
\end{figure}

\subsection{Poisson problem}
In this example, the Poisson problem (cf. Section \ref{sec: Poisson problem}) is solved for the temperature $u$ on two different two-dimensional domains. The first domain $\varOmega_\text{L}=\{(-1,1)\times(-1,1)\} \setminus \{(0,1)\times(0,1)\}$, referred to as the L-Shape, is characterized by a re-entrant corner with an opening angle of $\beta=90^{\circ}$ and a given exact solution
 \begin{equation}
 \bar{u}=r^{\frac{2}{3}} \sin \tfrac{2 \phi - \pi}{3}
 \label{eq:AnalSolHeat1}
\end{equation} in polar coordinates $(r,\phi)$.
The second domain $\varOmega_\text{S}=\{(-1,1)\times(-1,1)\}$, referred to as the slit domain, is characterized by a re-entrant corner with an opening angle of $\beta=0^{\circ}$ and a given exact solution
 \begin{equation}
 \bar{u}=r^{\frac{1}{2}} \sin \tfrac{\phi}{2}~.
 \label{eq:AnalSolHeat2}
\end{equation}     
Both boundary value problems are illustrated in Figure~\ref{fig:HeatConductionAufbau_1}. The boundary conditions are applied by setting $u=0$ at the Dirichlet boundary $\varGamma_\text{D}$ and the exact heat flux $g=\partial \bar{u} / \partial \nu_\text{N}$ at the Neumann boundary $\varGamma_\text{N}$. The L-Shape is modelled by a single $C^1$-continuous B-spline patch, while the slit domain is modelled by a single B-spline patch with $C^0$-continuous lines at the axis of symmetry of the domain as indicated by the dashed lines in Figure~\ref{fig:HeatConductionAufbau_1}.

In both problems, the geometry leads to a singularity of the solution at the re-entrant corner. In this case classical convergence theory does not hold, and the order of convergence with respect to the total number of degrees of freedom 
\begin{equation}
 k=- \tfrac{1}{2} \; \min \bigl( p,\tfrac{\pi}{2 \pi-\beta} \bigr)
\end{equation}
is governed by the angle $\beta$ of the re-entrant corner \cite{Yosibash2011}. For uniform $h$-refinement this leads to a convergence rate of 
$k=-1/3 \;\;\forall p$ for the L-shape and $k=-1/4 \;\;\forall p$ for the slit domain. 

The optimal order of convergence $k=-p/2$ can be recovered by local mesh refinement in the vicinity of the singularity. In the following, the adaptive finite element method (cf. Section \ref{sec: adaptive loop}) will be applied to solve the problem above with different refinement strategies. To select elements for refinement, the quantile marking (cf. Section~\ref{subsec: marking}) is used. The associated parameter $\alpha$ is adjusted for each refinement strategy, to achieve best possible convergence rates. 
\begin{figure}
 \centering
 \includegraphics[scale=0.80]{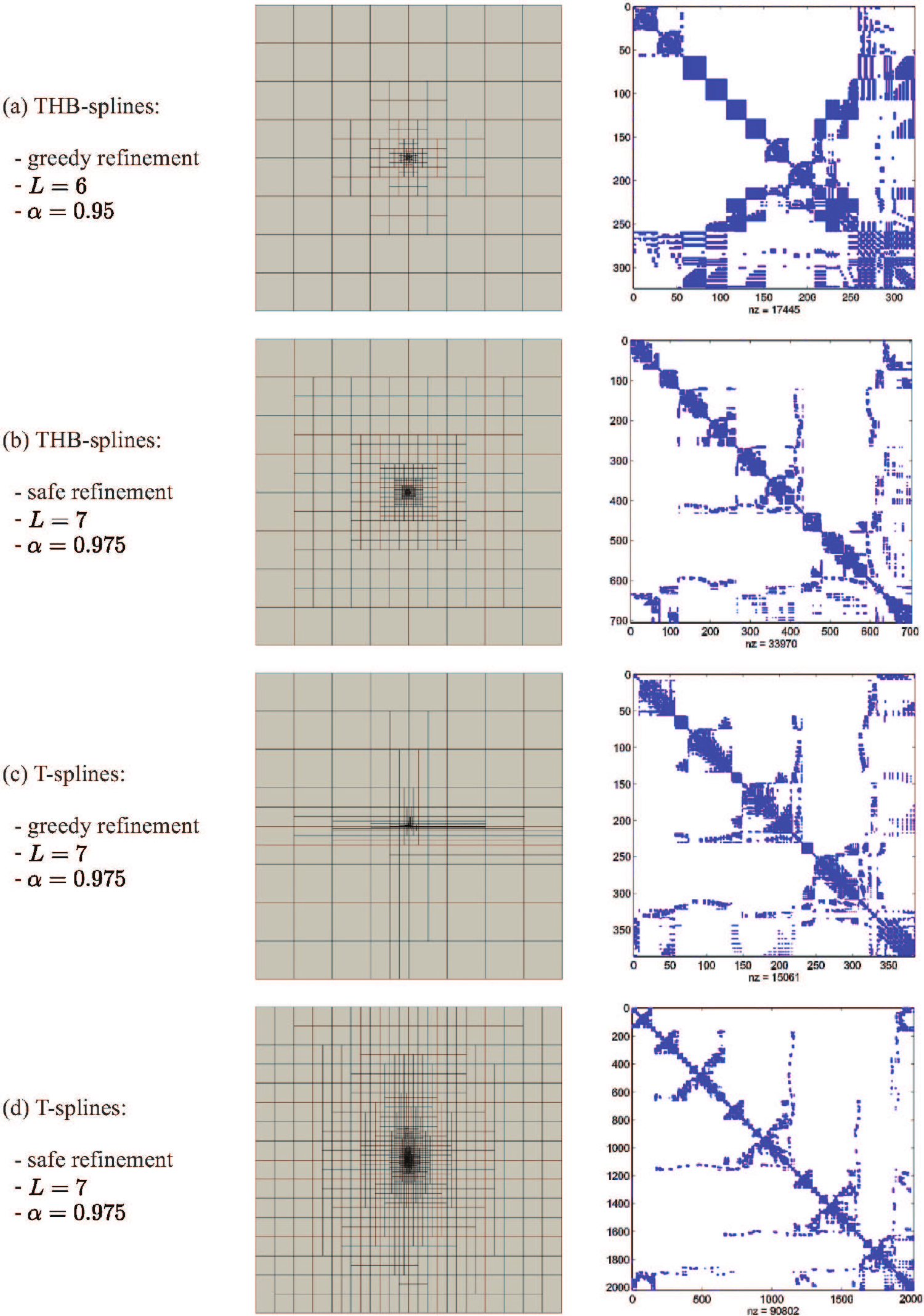}
 \caption{Slit domain: The marking parameters $\alpha$, the B\' ezier meshes and the sparsity patterns of the stiffness matrices after $L$ refinement steps for all \textbf{(a)}-\textbf{(d)} refinement strategies. The \much refinement strategies result in well graded meshes. Especially the \little T-spline refinement creates an unstructured mesh with badly shaped elements. Again, the \little THB-spline refinement creates the stiffness matrix with the highest density and interaction.}
 \label{fig:HeatConductionSD_1}
\end{figure}
\begin{figure}
 \centering
 \includegraphics[scale=0.70]{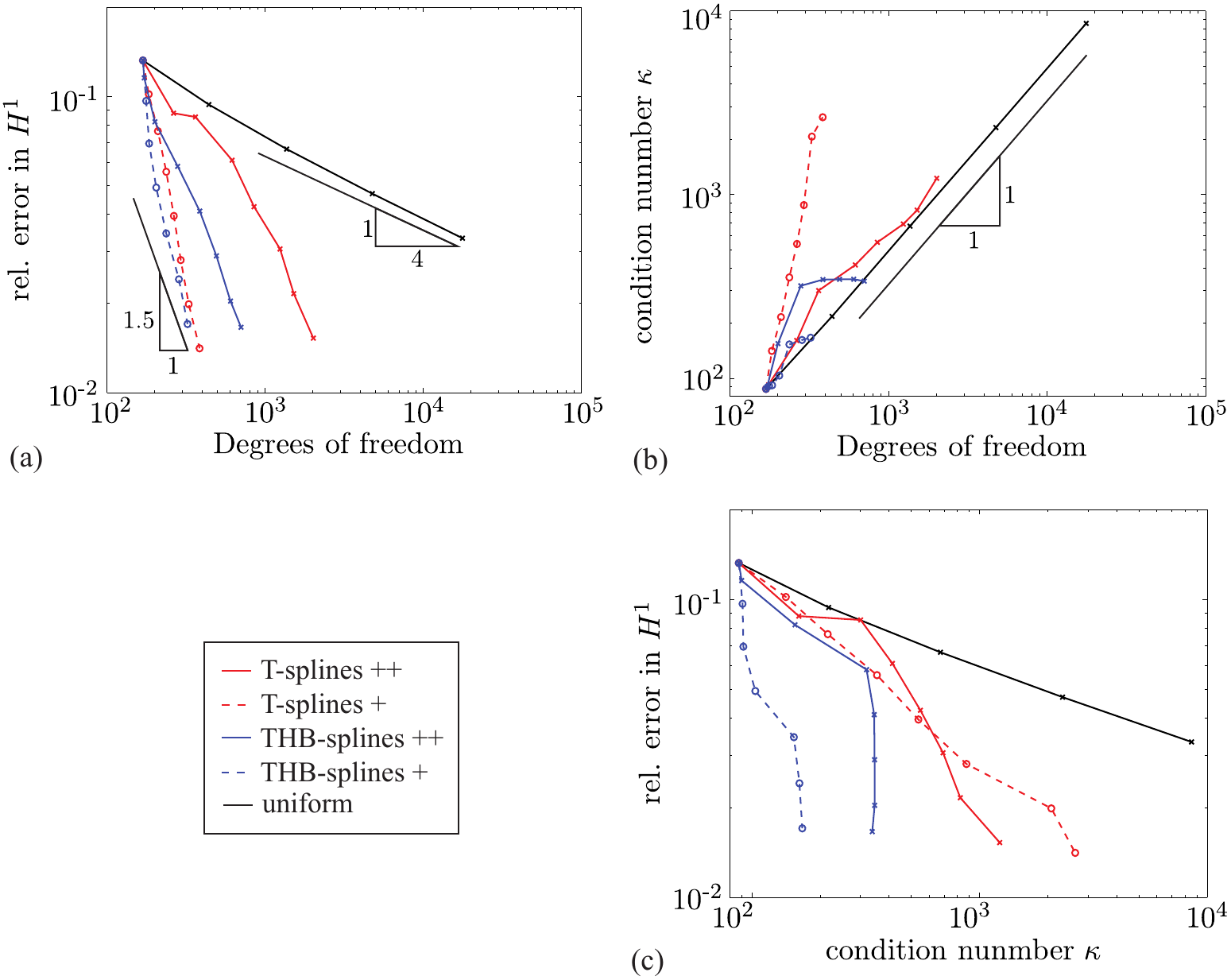}
 \caption{Slit domain: The convergence rates as well as the relations between the condition number of the stiffness matrix, the numerical error of the solution and the total number of degrees of freedom are illustrated. \textbf{(a)} - All refinement strategies converge with the expected convergence rate $k=1.5$ in the asymptotic range.}
 \label{fig:HeatConductionSD_2}
\end{figure}

\subsubsection{L-Shape}

The initial mesh of the L-shape problem consists of 16 elements. Figure~\ref{fig:HeatConductionLShape_1} shows the B\' ezier meshes after $L$ refinement steps, as well as the marking parameters $\alpha$. For the adaptive local refinement, the error in the $H^1$ norm is plotted over the total number of degrees of freedom (DOF) in Figure~\ref{fig:HeatConductionLShape_2}~(a). All refinement strategies recover the optimal order of convergence in the asymptotic range. Due to the coarse initial mesh, the \much refinements produce a greater amount of DOF in the pre-asymptotic range which is in particular observed for the \much T-spline refinement. As a result, the \much refinements are not as local as the \little refinements but create more smoothly graded meshes. To counteract the non-local refinements, the marking parameter for \much refinements is chosen higher.

Especially for the \little THB-spline refinement, the computed stiffness matrix has a higher density. For all other refinement strategies no clear tendency is visible in the sparsity patterns in Figure~\ref{fig:HeatConductionLShape_1}. 

The condition number is plotted over the DOF in Figure~\ref{fig:HeatConductionLShape_2}~(b). Due to the geometric map of the L-shape, a rate higher than one is reached for uniform refinement. Regarding the local refinement, results similar to the previous example are obtained. 
However, the differences between the \little and \much refinements are not as large as in the first experiment. 

As mentioned above, also the error of the numerical solution with respect to the condition number (cf. Figure~\ref{fig:HeatConductionLShape_2}~(c)) is of interest. It can be seen that for the same order of accuracy, all local refinement techniques produce smaller condition numbers compared to the uniform case. This means, that for local refinement, the error of the solution decreases faster per DOF than the condition number increases per DOF. This is an important result, because it illustrates that the negative influence of a locally refined mesh on the condition number does not predominate the benefits of local refinement regarding the error level. The refinement strategies compared among themselves show similar results.

\subsubsection{Slit domain}

The initial mesh of the slit domain consists of 64 elements. The B\' ezier meshes after $L$ refinement steps, as well as the marking parameters $\alpha$ are illustrated in Figure~\ref{fig:HeatConductionSD_1}. As expected, the meshes of the \much refinement routines propagate the refinement area but produce well graded meshes. On the other hand, the \little T-spline refinement leads to a mesh with little structure and badly shaped elements with aspect ratios up to 64. Concerning the sparsity patterns of the stiffness matrix, only the \little THB-spline refinement creates matrices with a higher density, due to the increased interaction between the basis functions. 

For the adaptive local refinement, the error in the $H^1$ norm is plotted over the total number of degrees of freedom (DOF) in Figure~\ref{fig:HeatConductionSD_2}~(a). It can be seen that the error of the \little refinement routines appear to converge with a higher rate in the pre-asymptotic range and later approach the theoretically predicted rate of $k=1.5$. The \much refinement routines have a minor convergence rate in the pre-asymptotic range, but then also converge with the theoretical rate of $k=1.5$. A reason for this behaviour can be found again in the relatively coarse initial mesh, which forces the \much T-spline refinement to refine almost the whole domain in the first refinement steps. As a result, the \much T-spline refinement requires six times more degrees of freedom than the \little T-spline refinement for the same error level.

The condition number is plotted over the DOF in Figure~\ref{fig:HeatConductionSD_2}~(b). Due to the badly shaped elements, the condition number for the \little T-spline refinement increases fastest. 
The THB-spline refinements instead seem to benefit from their hierarchical structure together with the absence of a deforming geometry mapping. At a certain stage of refinement, the condition number does not increase further. This behaviour has been also found in \cite{Johannessen2015} where HB-splines are compared against THB- and L-RB-splines.
In the context of hierarchical Finite Elements \cite{yserentant1986}, it is known and even proven that the condition number of the stiffness matrix scales with $O(\log(\texttt{DOF}))$ instead of $O(\texttt{DOF})$,  due to orthogonalities w.r.t.\ to the energy product between basis functions of different levels. In 1D, this leads to block-diagonal stiffness matrices; in higher dimensions, this effect is milder (see e.g.\ Figure~\ref{fig:HeatConductionSD_1}~(a)), but still yields good conditioning. It seems that (Truncated) Hierarchical B-splines share these benefits, however further investigation is needed in future.

Due to this effect, the \little THB-spline refinement performs best if the numerical error is plotted over the condition number (cf. Figure~\ref{fig:HeatConductionSD_2}~(c)). Since only a small amount of DOF is added during the refinement and due to the fact that the condition number grows slowly per DOF, an increased level of accuracy can be reached without increasing the condition number. But compared to the uniform refinement, also the T-spline refinements produce smaller condition numbers. 
 \begin{figure}
 \centering
 \includegraphics[scale=0.8]{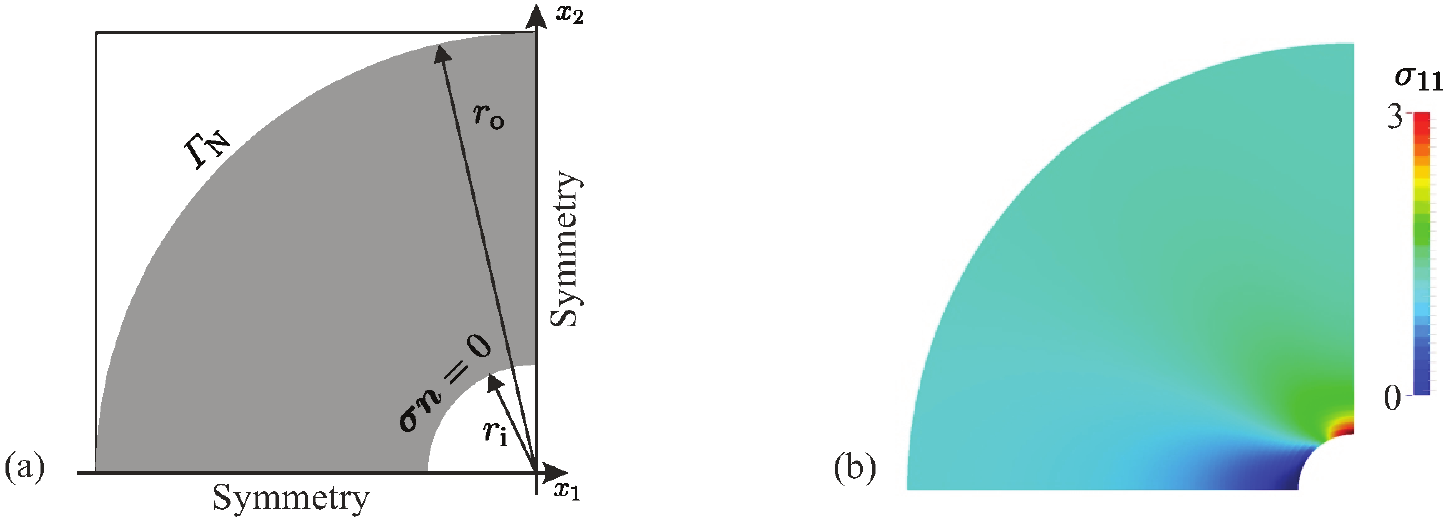}
 \caption{Infinite plate with a circular hole: \textbf{(a)} numerical analysis domain and boundary conditions, and \textbf{(b)} solution for $\sigma_{11}$.}
 \label{fig:PlateHoleR8_1}
\end{figure}  

\subsection{Linear elasticity}

As a third example, an infinite plate with a circular hole under uniaxial in-plane tension $\sigma_0$ according to Figure~\ref{fig:PlateHoleR8_1}~(a) is considered. 
The analytical solution is given by Timoshenko \cite{Timoshenko1970} in polar coordinates $(r,\phi)$
\begin{align}
 \bar{\sigma}_r & =\frac{\sigma_0}{2}\left[1-\frac{r_\text{i}^2}{r^2}+\left(1-4\frac{r_\text{i}^2}{r^2}+3\frac{r_\text{i}^4}{r^4}\right)\cos(2\varphi)\right]\\
 \bar{\sigma}_\varphi & =\frac{\sigma_0}{2}\left[1+\frac{r_\text{i}^2}{r^2}-\left(1+3\frac{r_\text{i}^4}{r^4}\right)\cos(2\varphi)\right]\\
 \bar{\sigma}_{r\varphi} & =\frac{\sigma_0}{2}\left(-1-2\frac{r_\text{i}^2}{r^2}+3\frac{r_\text{i}^4}{r^4}\right)\sin(2\varphi)
 \label{PlateHoleAnaSolu}
\end{align}
where $r_\text{i}=\SI{1}{\milli\meter}$ is the radius of the hole. A numerical solution is conveniently obtained on the quarter of an annulus with Dirichlet boundaries 
to enforce the symmetry conditions, and a Neumann boundary $\varGamma_\text{N}$ at the outer radius to enforce the exact normal stress.
The uniaxial tensile stress $\sigma_0=\SI{1}{\mega\pascal}$ is applied in the $x_1$-direction and material parameters $E=\SI{e5}{\pascal}$ and $\nu=0.3$ are used. 
The computational domain is modelled by a single $C^1$-continous NURBS patch with an outer radius $r_\text{o}=8$.  

\begin{figure}
 \centering
 \includegraphics[scale=0.80]{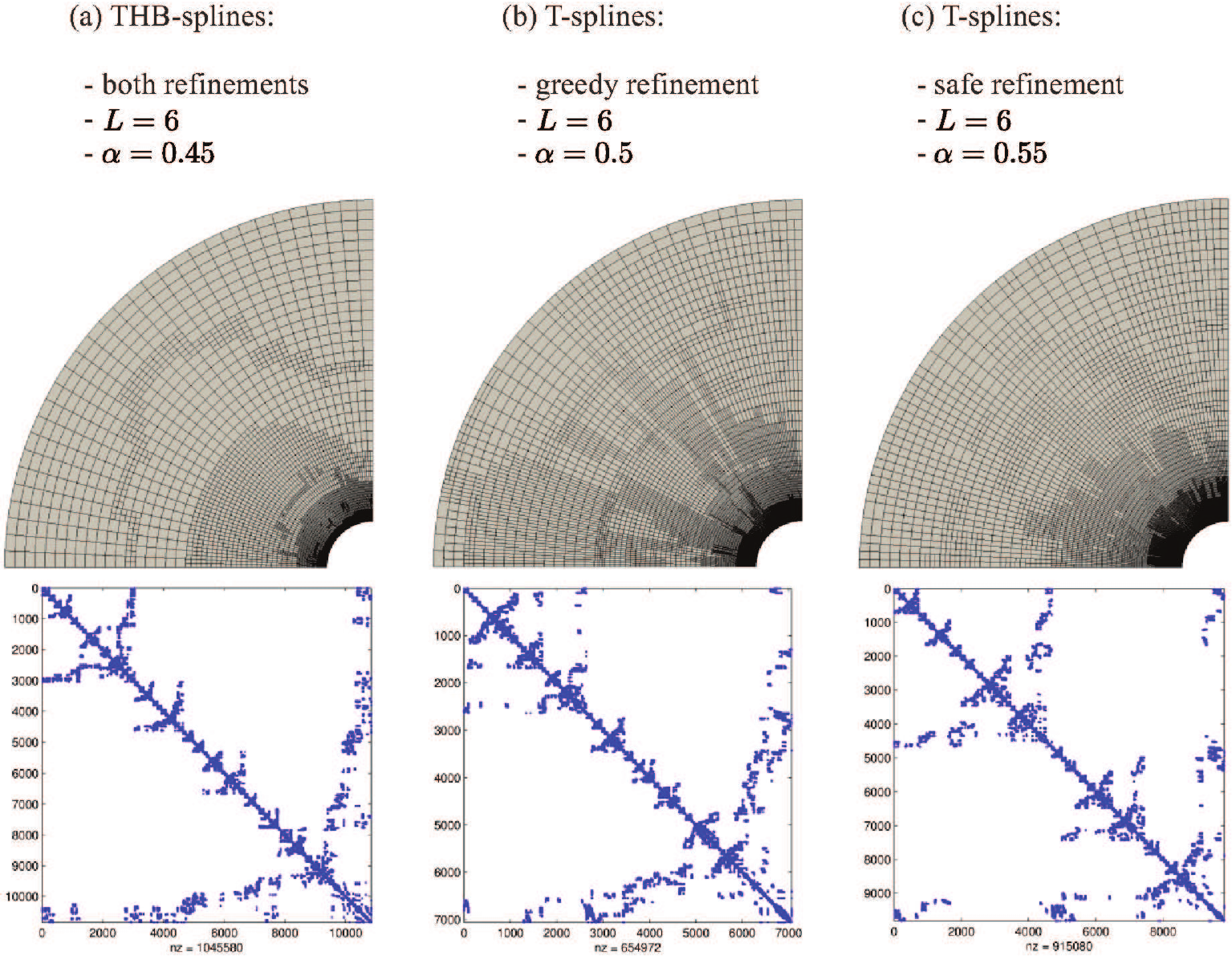}
 \caption{Infinite plate with a circular hole: The marking parameters $\alpha$, the B\' ezier meshes and the sparsity patterns of the stiffness matrices after $L$ refinement steps for all \textbf{(a)}-\textbf{(d)} refinement strategies. The \little and \much THB-spline refinement show an identical refinement behaviour. Neither in the B\' ezier meshes, nor in the sparsity patterns, clear differences between the refinement strategies are visible.}
 \label{fig:PlateHoleR8_2}
\end{figure}

The exact solution features a stress concentration at $(x,y)=(0,r_i)$ of $\sigma_{11}=3\sigma_0$ as illustrated in Figure~\ref{fig:PlateHoleR8_1}~(b). 
Due to the lack of a singularity, optimal convergence rates $k=-p/2$ can be obtained by uniform $h$-refinement. Local refinement does not improve this rate in the
asymptotic limit \cite{doerfel2010,Vuong2011}. There is however a benefit of the adaptive refinement which increases with the locality 
of the stress concentration. That is, if the outer radius $r_\text{o}$ is larger, the stress concentration is more localised in the computational domain, cf. Figure~ \ref{fig:PlateHoleR8_2}~(a), and an improved convergence rate can be achieved in the pre-asymptotic region. 

\begin{figure}
 \centering
 \includegraphics[scale=0.70]{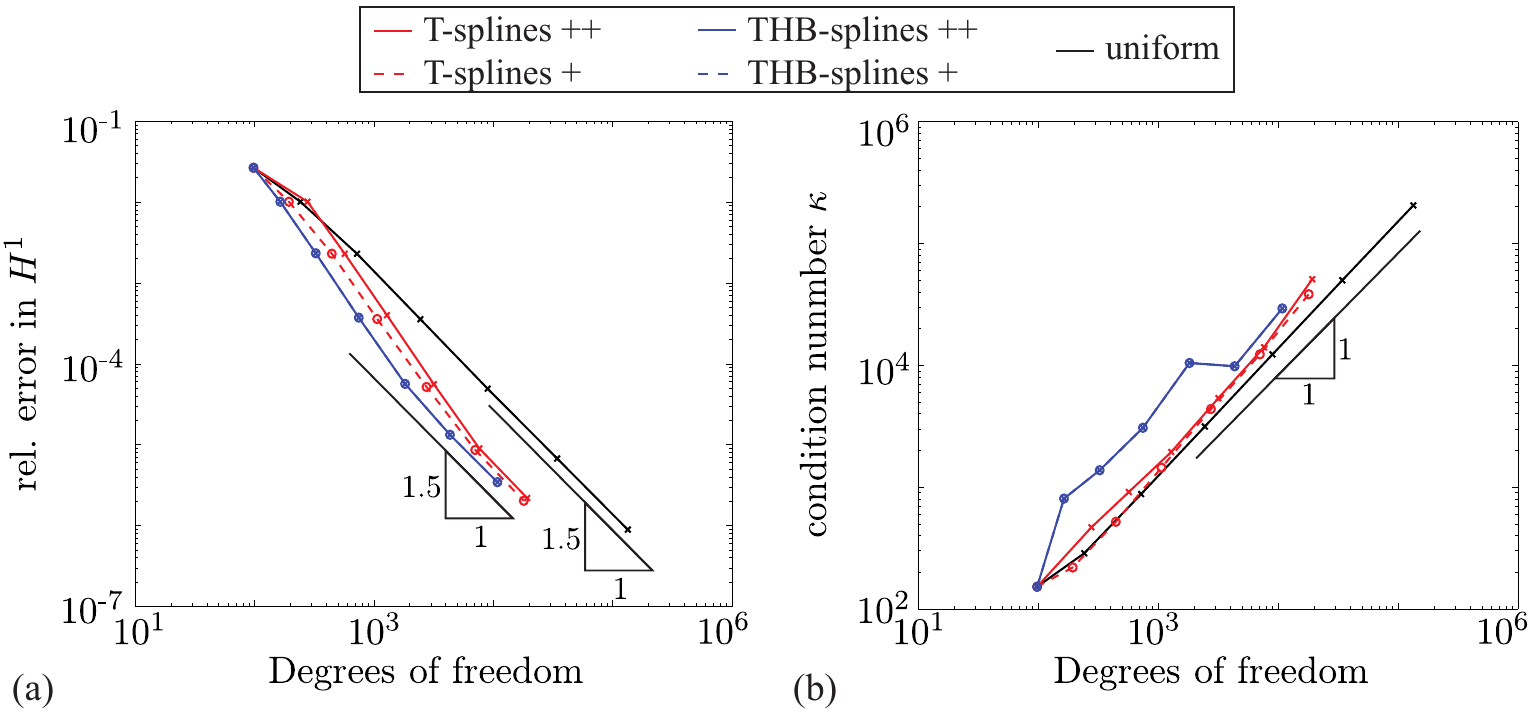}
 \caption{Infinite plate with a circular hole: The convergence rates as well as the condition number over the total number of degrees of freedom is plotted. \textbf{(a)} - The local refinement stategies improve the convergence rate in the pre-asymptotic regime, but reduce to $k=1.5$ in the asymptotic region.}
 \label{fig:PlateHoleR8_3}
\end{figure}

This improvement can be obtained for all refinement techniques by setting the marking parameter around $\alpha=0.5$ to generate a more extensive refinement. For this example the \little and \much THB-spline refinement produce same results. The meshes after $L$ refinement steps and the marking parameters $\alpha$ are illustrated in Figure~\ref{fig:PlateHoleR8_2}. All refinement techniques lead to similar meshes. As a result, also the sparsity patterns are quiet similar and do not show any tendency. If the condition number is plotted over DOF (cf. Figure~\ref{fig:PlateHoleR8_3}~(b)), no differences in the rate are visible between local and uniform refinement. In general it may be said that no numerical differences occur between the refinement techniques if the refinement area is extensive.  

\section{Conclusions}\label{sec: conclusion}
In this contribution, four different refinement techniques based on T- and THB-splines have been applied to the Adaptive Finite Element Method and compared regarding there theoretical and numerical properties. For this propose, four numerical examples have been studied.

In general, the successive use of an elementary refinement routine such as $\refthb$ or $\reftel$ causes uncontrolled function overlap and dense stiffness matrices in the case of THB-splines, or yields non-nested discrete spaces in the case of T-splines (which means for an isogeometric method that the geometry is not preserved), or does not even yield an actual basis due to linear dependencies (also for T-splines).

The refinement routine $\reftsp\p$ eliminates these major drawbacks from a practical approach, yielding an efficient and flexible refinement routine.
The procedure $\refthb\p$ is a practical approach to avoid the above-mentioned dense stiffness matrices, which is only partially achieved. On the other hand, it satisfies the same theoretical properties (namely linear complexity and bounded overlay) as the \much refinement routines, while such analysis for $\reftsp\p$ is currently not available.

The \much refinement routines $\refthb\pp$ and $\reftsp\pp$ have shown the expected optimal asymptotical behavior that has been predicted in theory, however they did not outperform the \little refinement routines in our experiments. 

Concerning the mesh grading and the numerical properties of the stiffness matrix, obvious differences increase with the locality of the problem. The refinement routines behave similar for extensive refinements, but differ the more local the refinement area is selected. For these local refinement areas, the \little refinement routines show a clear increase in the condition number per degree of freedom and especially for the refinement routine $\refthb\p$ dense stiffness matrices arise. Furthermore, the \little refinement routines lead to unstructured meshes around the refinement area. For the refinement routine $\reftsp\p$ this can lead to badly shaped elements with ever-growing aspect ratios (in our experiments, up to 64). 

From the implementation point of view, which is only a subjective view of the authors, the implementation effort increases from the refinement routine $\refthb\p$ to $\refthb\pp$ to $\reftsp\pp$ to $\reftsp\p$. The effort grows further for T-splines in general if a polynomial degree distinct from three is used, or for the \little T-spline refinement, if a generalisation to three dimension is desired.

%

%

\end{document}